%% file: main.tex
\documentclass[a4paper,11pt]{article}

\input{setup/praeamble}

\input{setup/colors}

\graphicspath{{./images/}}

\newtoggle{expensiveplots}

\toggletrue{expensiveplots}

\begin{document}

\author[,1]{Philipp L. Kinon \orcidlink{0000-0002-4128-5124} \footnote{Corresponding author: philipp.kinon@kit.edu}}
\author[2]{Tobias Thoma \orcidlink{0000-0002-6876-5662}}
\author[1]{Peter Betsch \orcidlink{0000-0002-0596-2503} }
\author[2]{Paul Kotyczka \orcidlink{0000-0002-6669-6368}}

\affil[1]{Institute of Mechanics, Karlsruhe Institute of Technology (KIT), Otto-Ammann-Platz 9, 76131 Karlsruhe, Germany}

\affil[2]{TUM School of Engineering and Design, Technical University of Munich (TUM), Boltzmannstr. 15, 85748 Garching, Germany}
\renewcommand\Authands{ and }

\input{includes/title}
\title{\mypapertitel{}}
\maketitle
\begin{abstract}
	\noindent Port-Hamiltonian (PH) systems provide a framework for modeling, analysis and control of complex dynamical systems, where the complexity might result from multi-physical couplings, non-trivial domains and diverse nonlinearities. A major benefit of the PH representation is the explicit formulation of power interfaces, so-called ports, which allow for a power-preserving interconnection of subsystems to compose flexible multibody systems in a modular way. In this work, we present a PH representation of geometrically exact strings with nonlinear material behaviour. Furthermore, using structure-preserving discretization techniques a corresponding finite-dimensional PH state space model is developed. Applying mixed finite elements, the semi-discrete model retains the PH structure and the ports (pairs of velocities and forces) on the discrete level. Moreover, discrete derivatives are used in order to obtain an energy-consistent time-stepping method. The numerical properties of the newly devised model are investigated in a representative example. The developed PH state space model can be used for structure-preserving simulation and model order reduction as well as feedforward and feedback control design.
\end{abstract}


\input{includes/body}

\FloatBarrier


\newenvironment{authcontrib}[1]{%
	\subsection*{\textnormal{\textbf{Author Contributions}}}%
	\noindent #1}%
{}%
\newenvironment{acks}[1]{%
	\subsection*{\textnormal{\textbf{Acknowledgements}}}%
	\noindent #1}%
{}%
\newenvironment{funding}[1]{%
	\subsection*{\textnormal{\textbf{Funding}}}%
	\noindent #1}%
{}%
\newenvironment{dci}[1]{%
	\subsection*{\textnormal{\textbf{Declaration of conflicting interests}}}%
	\noindent #1}%
{}%
\newenvironment{code}[1]{%
	\subsection*{\textnormal{\textbf{Code}}}%
	\noindent #1}%
{}%


%
\begin{acks}
	The financial support for this work by the DFG (German Research Foundation) – project number 388118188 – is gratefully acknowledged.
\end{acks}
\begin{dci}
	The authors declare no conflict of interest.
\end{dci}
\begin{code}
	The source code used for the finite element computations is implemented in Matlab and can be obtained on reasonable request.
\end{code}

\addcontentsline{toc}{section}{References}
\bibliographystyle{unsrtnat}
\bibliography{bib}

\end{document}

%% file: setup/praeamble.tex
\usepackage{moreverb,url}
\usepackage[
	colorlinks,
	bookmarksopen,
	bookmarksnumbered,
	citecolor=blue,
	urlcolor=blue,
	linkcolor=blue,
]{hyperref}
\usepackage[font=small,labelfont=bf]{caption}
\usepackage{subcaption}
\usepackage{csquotes}
\usepackage{mathtools}
\usepackage{xcolor}
\usepackage{amsfonts}
\usepackage{amsmath}
\DeclareMathSymbol{\shms}{\mathbin}{AMSa}{"39}  
\usepackage{etoolbox}	
\usepackage{adjustbox}
\usepackage{placeins}
\usepackage[title]{appendix}
\usepackage{float}
\usepackage{siunitx}
\usepackage{breqn}

\usepackage{authblk}		
\usepackage{orcidlink}

\usepackage{geometry}
\usepackage[numbers,sort]{natbib}

\renewcommand{\d}{{\,\rm  d}}

\newcommand\BibTeX{{\rmfamily B\kern-.05em \textsc{i\kern-.025em b}\kern-.08em
			T\kern-.1667em\lower.7ex\hbox{E}\kern-.125emX}}

\usepackage{accents}
\newlength{\dhatheight}

\newcommand{\transp}{^\mathrm{T}}

\renewcommand{\d}[1][]{\ensuremath{\,\mathrm{d}#1}}
\renewcommand{\phi}{\varphi}
\newcommand{\D}[1][]{\ensuremath{\mathrm{D}#1}}
\usepackage{color}
\usepackage{transparent}
\newcommand{\npe}{_{n+1}}
\newcommand{\npeh}{_{n+1/2}}
\newcommand{\n}{_{n}}
\newcommand{\h}{^\mathrm{h}}
\newcommand{\PDG}[1][]{\ensuremath{\bar{\partial}#1}}
\newcommand{\DG}[1][]{\ensuremath{\bar{\mathrm{D}}#1}}
\newcommand{\eqcomma}{\ensuremath{\text{,}}}
\newcommand{\eqdot}{\ensuremath{\text{.}}}
\renewcommand{\vec}[1]{\ensuremath{\mbox{$ #1 $}}}
\usepackage{bm}
\renewcommand{\epsilon}{\varepsilon}
\newsavebox{\topleftimage}
\newsavebox{\bottomleftimage}
\newsavebox{\bigimage}
\usepackage{mwe}
\usepackage[ruled]{algorithm2e}

\usepackage{pgfplots}
\pgfplotsset{compat = newest}
\usepackage{pgfplotstable}
\usepackage{tikz}
\usepackage{tikz-3dplot} 
\usetikzlibrary{arrows.meta,decorations.markings}
\usepgfplotslibrary{fillbetween}
\pgfkeys{/pgf/number format/.cd,fixed,precision=4}

\newlength\figH
\newlength\figW
\usepgfplotslibrary{patchplots}

\tikzset{myarrowhead/.style={decoration={markings,mark=at position 1 with %
{\arrow[scale=1.,>={Triangle[length=8pt, width=7pt]}]{>}}},postaction={decorate}}}

\DeclareUnicodeCharacter{2212}{−}
\usepgfplotslibrary{groupplots,dateplot}
\usepgfplotslibrary{polar}
\usetikzlibrary{positioning}
\usetikzlibrary{patterns,shapes.arrows}
\pgfplotsset{compat=newest}

\usetikzlibrary{external}

\pgfplotsset{every tick label/.append style={font=\scriptsize}}



%% file: setup/colors.tex
\definecolor{color1}{RGB}{230, 159, 0}
\definecolor{color2}{RGB}{86, 180, 233}
\definecolor{color3}{RGB}{204, 121, 167}
\definecolor{color4}{RGB}{0, 158, 115}
\definecolor{color5}{RGB}{0, 114, 178}
\definecolor{color6}{RGB}{213, 94, 0}
\definecolor{color7}{RGB}{240, 228, 66}
\definecolor{colorblack}{RGB}{0, 0, 0}          

\definecolor{cgrey}{RGB}{128, 128, 128}
\definecolor{cdarkgrey}{RGB}{98, 98, 98}

\colorlet{corange}{color1}
\colorlet{ccyan}{color2}
\colorlet{cviolet}{color3}
\colorlet{cgreen}{color4}
\colorlet{cblue}{color5}
\colorlet{cred}{color6}
\colorlet{cyellow}{color7}

\colorlet{ctop}{color1}
\colorlet{ccenter}{color2}
\colorlet{cbottom}{color3}

\colorlet{cpath1}{color1}
\colorlet{cpath2}{color2}
\colorlet{cpath3}{color3}
\colorlet{cpath4}{color4}
\colorlet{cpath5}{color7}

\colorlet{cpathneutral}{cgrey}

\colorlet{cpathm90}{color1}
\colorlet{cpathm45}{color2}
\colorlet{cpath0}{color3}
\colorlet{cpath45}{color4}
\colorlet{cpath90}{color7}

\colorlet{cmtoa}{cred}
\colorlet{cHS01}{cgreen}
\colorlet{cHS02}{color4}
\colorlet{cHS03}{cblue}
\colorlet{cmtlinear}{cyellow}
\colorlet{cmtlinearcompliance}{cviolet}

%% file: includes/title.tex
\newcommand{\mypapertitel}{%
	Port-Hamiltonian formulation and structure-preserving discretization of hyperelastic strings\footnote{submitted as a proceeding to the ECCOMAS Thematic Conference on Multibody Dynamics 2023}
}

%% file: includes/body.tex
\section{Introduction}
The class of port-Hamiltonian (PH) systems has become increasingly important across multiple research fields dealing with modeling and control of dynamical systems. PH systems provide a framework for the analysis of complex dynamical systems \cite{duindam_modeling_2009,schaft_port-hamiltonian_2014}, where the complexity might result from multi-physical couplings, non-trivial domains and diverse nonlinearities. A major benefit of the PH representation is the explicit formulation of power interfaces, so-called ports, which allow for an intrinsically power-preserving interconnection of subsystems. In this way, the modular composition of a model can be facilitated.

Since many technical (sub-)systems are modeled by partial differential equations (PDE), the theory of infinite-dimensional PH systems (see, e.g., \cite{van2002hamiltonian}) has been extended in recent years \cite{jacob2012linear,rashad2020twenty}. Among examples from various physical disciplines, PH formulations have also been proposed in structural mechanics, e.g., for flexible multibody dynamics \cite{brugnoli2021port}. However, many contributions focus on linear theories and linear material behaviour, e.g., \cite{warsewa2021port}.

The structural elements of strings (see, e.g., \cite{antman_nonlinear_2005}) are particularly interesting and are widely used in control and modeling of multibody systems (see, e.g., \cite{strohle2022simultaneous}). Strings can be found as a sub-module in a large variety of systems such as cranes, cable robots or satellite systems. Accordingly, they are interconnected with their environment, which can be described beneficially within the PH framework. To the best of our knowledge, the PH representations of strings have yet been restricted to linear stress-strain relations \cite{thoma2022port}. We try to fill this gap with the present contribution.

We propose a nonlinear port-Hamiltonian model for geometrically exact strings with hyperelastic constitutive relations. We thus enhance the previous work \cite{thoma2022port} with respect to material nonlinearities and use strain and stress measures inspired from nonlinear continuum mechanics. The formulation reveals power-conjugated interfaces, so-called ports, on the boundary of the strings and facilitates the energy-consistent coupling in multibody systems. The resulting formulation can be linked to the Hu-Washizu principle from the theory of elasticity (see, e.g., \cite{betsch_energy_2016}).

With respect to the numerical discretization, the mixed finite element method has been employed to establish an approximate model under preservation of the PH structure (as, e.g., in \cite{cardoso2016piezoelectric, cardoso2021partitioned}). A core characteristic of structure-preserving discretization of PH systems is retaining the ports (e.g., pairs of velocities and forces) with their causality (which refers to the definition of boundary input variables in the system theoretic sense) on the discrete level. We moreover apply the discrete gradient in the sense of Gonzalez \cite{gonzalez_time_1996} to achieve an energy-consistent discretization in time, which inherits the passivity property of the underlying continuous model.

The work is structured as follows. In Sec.~\ref{sec_string} the governing equations of nonlinear strings are recapitulated. We propose a mixed formulation in the port-Hamiltonian framework and analyze the continuous power balance. Sec.~\ref{sec_discrete} shows our chosen discretization methods, both in space (mixed finite elements) and time (discrete gradients) in order to obtain a discrete PH model for nonlinear strings. The resulting numerical properties are studied in a representative example in Sec.~\ref{sec_numerics} before concluding our work in Sec.~\ref{sec_conclusion} along with some outlook for future research directions.

\section{Port-Hamiltonian string representation}\label{sec_string}
In this section, the fundamental governing equations for geometrically exact strings are refreshed. Furthermore, a new port-Hamiltonian formulation for these equations is introduced.

\subsection{Governing equations}

\begin{figure}[tb]
    \centering
    \def\svgwidth{0.6\textwidth}
    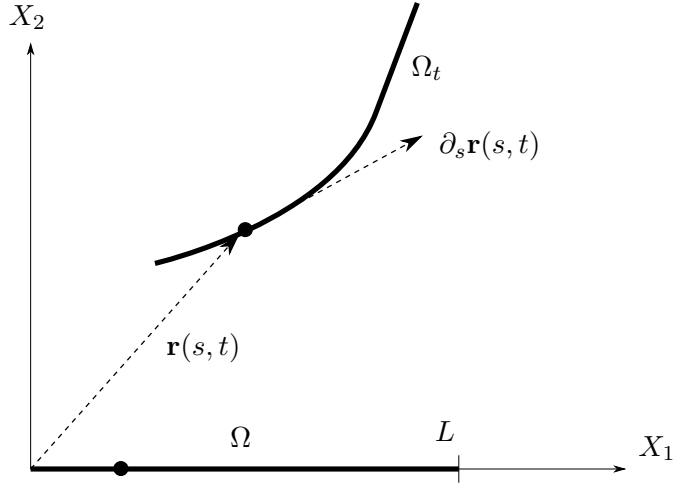
    \caption{Material and spatial string configurations on $\mathbb{R}^2$.}
    \label{fig:configuration}
\end{figure}

We consider a one-dimensional undeformed (material) string configuration $\Omega \subset \mathbb{R}^d$ of length $L\in\mathbb{R}^+$ and its current (spatial) configuration $\Omega_t\subset\mathbb{R}^d$, where $d \in \{1, 2, 3\}$ is the spatial dimension. Correspondingly, the position vector $\vec{r}: S \times T \rightarrow \Omega_t$ is introduced, depending on the material arc-length coordinate $s \in S = [0,L] \subset \mathbb{R}_0^+$ and on time $t \in T \subset \mathbb{R}_0^+$, see Fig.~\ref{fig:configuration} for an illustration on $\mathbb{R}^2$. The balance of linear momentum of the string (see, e.g., \cite{strohle2022simultaneous}) in material coordinates is given by
\begin{align}
    \label{balance_momentum}
    \rho A \ddot{\vec{r}}(s,t) = \partial_s \vec{n}(s,t) + \vec{b}
\end{align}
and includes the constant density per unit length $\rho A \in\mathbb{R}_0^+$, the constant body forces per unit length $\vec{b} \in \mathbb{R}^d$. Moreover, $\partial_s(\square)$ represents the partial derivative with respect to the material coordinate and the contact force per unit length $\vec{n}(s,t) : S \times T \rightarrow  \mathbb{R}^d$ is introduced. Due to the kinematic assumption (i.e., neglect bending stiffness), the contact force $\vec{n}$ is oriented tangentially such that
\begin{align}
    \vec{n}(s,t) = N(s,t) \frac{\partial_s \vec{r}(s,t)}{||\partial_s \vec{r}(s,t||} , \label{normal_force}
\end{align}
where the tension $N: S \times T \rightarrow \mathbb{R}$ follows from a constitutive relation as presented in Sec.~\ref{sec_2_2}. We consider hyperelastic materials with a stored energy density $W : \mathbb{R} \rightarrow \mathbb{R}$ depending on the strain type quantity
\begin{align}
    C=\partial_s \vec{r} \cdot \partial_s \vec{r}, \label{rightCG}
\end{align}
which can be interpreted in analogy to the right Cauchy-Green strain tensor from elasticity. Accordingly, the conjugated stress quantity (in analogy to the second Piola-Kirchhoff stress tensor) can be obtained via
\begin{align}
    S := 2 \D W(C) \label{const_relation} ,
\end{align}
where $\D (\square)$ denotes the gradient operator. We refer to the textbook by Antman \cite{antman_nonlinear_2005} for detailed analyses of strings. Further kinematic assumptions for arriving at the balance of linear momentum \eqref{balance_momentum} are well summarized in \cite{thoma2022port}. From now on, we drop both spatial and temporal arguments for the sake of brevity and only use them where the explicit mention is helpful.

\subsection{Constitutive modeling}\label{sec_2_2}

As shown in \cite{strohle2022simultaneous} a stored energy density function for strings with hyperelastic material is given by
\begin{align}
    U(\nu) = \frac{EA}{4} (\nu^2 - 2 \ln(\nu) -1),
\end{align}
where the stretch $  \nu = || \partial_s \vec{r} || = \sqrt{C}$ is used and $EA$ denotes the constant axial stiffness. The energy density describes the stored energy per unit length due to elastic deformation. Consequently, the tension follows directly by differentiation, i.e.,
\begin{align}
    N(\nu) = \D U(\nu).
\end{align}
In the present work, however, the strain quantity \eqref{rightCG} shall be used. Stating $W(C) = U(\nu(C))$, we obtain
\begin{align} \label{const_relation_hyper}
    W(C) = \frac{EA}{4} (C - 2 \ln(\sqrt{C}) -1) = \frac{EA}{4} (C - \ln(C) -1) \eqdot
\end{align}
This yields the stress
\begin{align}
    S = 2 \D W(C) = \frac{EA}{2} (1 - C^{-1}) =\frac{EA}{2C} (C -1) \eqdot
\end{align}
It is moreover crucial to find a relation between the tension $N$ introduced in \eqref{normal_force} and the corresponding stress $S$.
To this end, relation \eqref{normal_force} can be rewritten such that
\begin{align}
    \vec{n} = N(\nu) \frac{\partial_s \vec{r}}{|| \partial_s \vec{r} ||} = \D U(\nu) \frac{\partial_s \vec{r}}{\sqrt{C}} = \D W(C) \D C(\nu) \frac{1}{\sqrt{C}} \partial_s \vec{r}
\end{align}
and eventually
\begin{align}
    \vec{n}(s,t) = S \partial_s \vec{r} .
\end{align}

\, \\
\noindent \textbf{Remark 1 (Linear elasticity):} The choice of a strain energy density function
\begin{align}
    W = \breve{W}(C) = \frac{1}{2} EA (\sqrt{C} -1)^2 = \frac{1}{2} EA \epsilon^2 \label{linear_elastic},
\end{align}
with $\epsilon := \sqrt{\partial_s \vec{r} \cdot \partial_s \vec{r}} -1 = \sqrt{C} -1$, corresponds to the previous work \cite{thoma2022port}. Hence, the linear elastic case is covered by the present approach.

\, \\
\noindent \textbf{Remark 2 (St.~Venant-Kirchhoff material):} A suitable representation of a linear stress-strain relation in the case of finite displacements (i.e., St. Venant-Kirchhoff material) is given by
\begin{align}
    W = \tilde{W}(C) = \frac{1}{8} EA (C - 1)^2 \label{st_venant}
\end{align}
and thus the stress reads
\begin{align}
    S = 2 \D \tilde{W}(C) = \frac{1}{2} EA (C-1) \eqdot
\end{align}
The stored energy density function \eqref{st_venant} can be deduced by applying the 1D string kinematic assumptions to the 3D material law of St. Venant-Kirchhoff.

\subsection{Port-Hamiltonian formulation}

With the above assumptions, the total energy of a geometrically exact string is described by its Hamiltonian
\begin{align}
    \label{hamiltonian}
    H(\vec{x})=H(\vec{r}, \vec{v},C) = T(\vec{v}) + V^{\mathrm{int}}(C) + V^{\mathrm{ext}}(\vec{r}) = \int_{0}^{L} \left( \frac{1}{2} \rho A\vec{v} \cdot \vec{v} + W(C) - \vec{r}\cdot \vec{b} \right) \d s .
\end{align}
Therein, $T$ denotes the kinetic energy and $V^{\mathrm{int}}$ and $V^{\mathrm{ext}}$ are the internal and external potential energy, respectively. In this model the state variables $\vec{x}$ contain positions $\vec{r}$, velocities $\vec{v} \in\mathbb{R}^3$ and the strain \eqref{rightCG}. The co-state variables (efforts) emerge as the variational derivatives
\begin{subequations}
    \begin{align}
        \frac{\delta H}{\delta \vec{r}} & = -\vec{b} \eqcomma              \\
        \frac{\delta H}{\delta \vec{v}} & = \rho A\vec{v} \eqcomma         \\
        \frac{\delta H}{\delta C}       & = \D W(C) = \frac{1}{2} S \eqdot
        \label{efforts}
    \end{align}
\end{subequations}
Taking the balance of linear momentum \eqref{balance_momentum}, the kinematic relation between the velocity $\vec{v}$ and the position $\vec{r}$, and the strain rate derived from \eqref{rightCG},
the state equations governing the motion of the string can be written as partial differential equations of first order in time
\begin{align}
    \label{eqn_of_motion}
    \begin{bmatrix}
        \vec{I} & 0             & 0 \\
        0       & \rho A\vec{I} & 0 \\
        0       & 0             & 1
    \end{bmatrix}
    \begin{bmatrix}
        \dot{\vec{r}} \\ \dot{\vec{v}} \\ \dot{C}
    \end{bmatrix} & = \begin{bmatrix}
                          0        &  & \vec{I}                                       &  & 0                                             \\
                          -\vec{I} &  & 0                                             &  & 2 \partial_s ( \square \, \partial_s \vec{r}) \\
                          0        &  & 2 \partial_s \vec{r} \cdot \partial_s \square &  & 0
                      \end{bmatrix}
    \begin{bmatrix}
        -\vec{b} \\
        \vec{v}  \\
        \frac{1}{2} S
    \end{bmatrix} \eqdot
\end{align}
These equations can be recast in the compact form
\begin{equation} \label{eq:PH-continuous-Mehrmann-Morandin2019}
    \mathcal{E}\dot{\vec{x}} = \mathcal{J}\vec{z} \ ,
    \quad\mbox{where}\quad
    \mathcal{E}\transp\vec{z} = \delta_{\vec{x}} H \ ,
\end{equation}
that lies at the heart of the PH formulation in \cite{beattie2018linear}, see also \cite{mehrmann19, brugnoli2021port}. Note that $\mathcal{E}\transp = \mathcal{E}$ and that the operator $\mathcal{J}$ is formally skew-adjoint ($\mathcal{J}^* = -\mathcal{J}$), which is proven by $\int_0^L \vec{z} \cdot \mathcal{J}\vec{z} \d s = 0 $ under zero boundary conditions, see \eqref{power_balance}.
In order to close the PH system representation, boundary conditions defining the system input $\vec{u}(t)$ and the corresponding collocated output $\vec{y}(t)$ are required. In the following, we consider a string system which is purely subject to Neumann boundary conditions (i.e., prescribed boundary contact forces), such that
\begin{align}
    \vec{u}(t) =
    \begin{bmatrix}
        -\vec{n}(0,t) \\ \vec{n}(L,t)
    \end{bmatrix}.
\end{align}
Following the energy balance
\begin{align}
    \dot{H} & = \int_{0}^{L} \delta_{\vec{x}} H \cdot \dot{\vec{x}} \d s = \int_{0}^{L} \mathcal{E}\transp\vec{z} \cdot \dot{\vec{x}} \d s =
    \int_{0}^{L} \vec{z} \cdot \mathcal{E}\dot{\vec{x}} \d s =
    \int_{0}^{L} \vec{z} \cdot \mathcal{J}\vec{z} \d s
    \label{power_balance}                                                                                                                                                                                                                                                   \\
            & = \int_{0}^{L} \left( \vec{v} \cdot \partial_s (S \partial_s \vec{r}) + S \partial_s \vec{r} \cdot \partial_s \vec{v} \right) \d s = \int_{0}^{L} \partial_s (\vec{v} \cdot S \partial_s \vec{r}) \d s = \left[ \vec{n} \cdot \vec{v} \right]_0^L      \notag \\
            & = \vec{u} \cdot \vec{y} \notag
\end{align}
immediately leads to the power conjugated output, given by the boundary velocities,
\begin{align}
    \vec{y}(t) = \begin{bmatrix}
                     \vec{v}(0,t) \\
                     \vec{v}(L,t)
                 \end{bmatrix} .
\end{align}
Hence, relation \eqref{power_balance} demonstrates the passivity and losslessness of the system, since the total change of energy is given by the power transmitted through the boundary $\partial\Omega=\{0,L\}$. This includes energy-conservation if $\vec{u} = \vec{0}$. In order to define the initial value problem, initial conditions
\begin{subequations}
    \begin{align}
        \vec{r}(s,0) & = \vec{r}_0(s), \\
        \vec{v}(s,0) & = \vec{v}_0(s), \\
        C(s,0)       & = C_0(s)
    \end{align}
\end{subequations}
are required, where $\vec{r}_0$, $\vec{v}_0$ and $C_0$ have to be prescribed, with $\vec{r}_0(s)$ and $C_0(s)$ satisfying relation \eqref{rightCG}.

\, \\
\noindent \textbf{Remark 3 (Mixed boundary conditions):} In case of mixed boundary conditions, the boundary is split into a Dirichlet and Neumann part, $\partial\Omega = \partial\Omega_N\cup\partial\Omega_D$, with $\partial\Omega_N\cap\partial\Omega_D = \emptyset$. Therefore, velocity inputs (and corresponding positions) are enforced as Dirichlet and force inputs as Neumann boundary conditions. Accordingly, the power conjugated output corresponds to the reaction force satisfying the Dirichlet condition or the velocities at the Neumann boundary.

\, \\
Towards the end of this section the weak form pertaining to \eqref{eqn_of_motion} is deduced. Standard procedures lead to the weak form, such that
\begin{subequations}
    \label{weak_form}
    \begin{align}
        \int_{0}^{L} \delta \vec{b} \cdot (\dot{\vec{r}} -\vec{v}) \d s                                                                                                                                                & = 0                     \\
        \int_{0}^{L} \left( \delta \vec{v} \cdot (\rho A \dot{\vec{v}} - \vec{b}) + \partial_s (\delta \vec{v}) \cdot S \partial_s \vec{r} \right) \d s - \left[ S \partial_s \vec{r} \cdot \delta \vec{v} \right]_0^L & = 0 \label{weak_form_2} \\
        \int_{0}^{L} \delta S \, ( \dot{C} - 2 \partial_s \vec{r} \cdot  \partial_s \vec{v}) \d s                                                                                                                      & = 0                     \\
        \int_{0}^{L} \delta C \, ( S - 2 \D \hat{W}(C)) \d s                                                                                                                                                           & = 0
    \end{align}
\end{subequations}
must hold for test functions $\delta \vec{b}, \delta \vec{v}, \delta S$ and $\delta C$ from appropriate spaces. Note that the constitutive relation \eqref{const_relation} has been appended in weak form to close the set of equations. This will be necessary to retain the PH structure in the space-discrete setting later on. Moreover, in the second equation, integration by parts has been used.
In the case of static problems, the related 3-field formulation in $(\vec{r},C,S)$ can be linked to the Hu-Washizu principle from the theory of elasticity (see \cite{betsch_energy_2016} and the inflated form in \cite{betsch_mixed_2018}).


\section{Structure-preserving discretization}
\label{sec_discrete}
In the sequel, the weak form \eqref{weak_form} is discretized by means of a mixed Galerkin finite element approach to obtain a semidiscrete state space model that fits again into the PH framework. Thereafter, a suitable time discretization is applied to obtain a structure-preserving scheme.

\subsection{Spatial discretization}
For the discretisation in space we divide the string domain into $n_{el}$ finite elements with respective domain $\Omega_e$ such that $\Omega = \cup_{e=1}^{n_{el}} \Omega_e$. Correspodingly, local ansatz functions are chosen and the isoparametric concept applies. In particular, we choose $C^0$-continuous, linear Lagrangian shapefunctions $\vec{\Phi}$ and elementwise constant, discontinuous ansatz functions $\vec{\Psi}$ such that the approximations are given by
\begin{subequations}
    \label{fe_approx}
    \begin{align}
        \vec{r}\h(s,t) & = \vec{\Phi}(s) \hat{\vec{r}}(t), & \delta \vec{b}\h(s) & = \vec{\Phi}(s) \delta \hat{\vec{b}}, \\
        \vec{v}\h(s,t) & = \vec{\Phi}(s) \hat{\vec{v}}(t), & \delta \vec{v}\h(s) & = \vec{\Phi}(s) \delta \hat{\vec{v}}, \\
        C\h(s,t)       & = \vec{\Psi}(s) \hat{\vec{C}}(t), & \delta C\h(s)       & = \vec{\Psi}(s) \delta \hat{\vec{C}}, \\
        S\h(s,t)       & = \vec{\Psi}(s) \hat{\vec{S}}(t), & \delta S\h(s)       & = \vec{\Psi}(s) \delta \hat{\vec{S}},
    \end{align}
\end{subequations}
where $\hat{(\square)}$ are arrays containing the nodal vectors of the respective quantity. Inserting the approximations \eqref{fe_approx} into the weak form \eqref{weak_form} yields the semi-discrete equations of motion
\begin{subequations} \label{eq:semi-discrete}
    \begin{align}
        \dot{\hat{\vec{r}}}              & = \hat{\vec{v}}  \eqcomma                                                            \\
        \vec{M}_\rho \dot{\hat{\vec{v}}} & = \vec{F}_b - \vec{K}(\hat{\vec{r}}) \hat{\vec{S}} +  \vec{F}_{\mathrm{N}}  \eqcomma \\
        \vec{M}_S \dot{\hat{\vec{C}}}    & = 2  \vec{K}(\hat{\vec{r}})\transp \hat{\vec{v}}   \eqcomma                          \\
        \vec{M}_S \hat{\vec{S}}          & = \int_{0}^{L} \vec{\Psi} \transp 2 \D \hat{W}(C\h) \d s \label{discrete_const1}
    \end{align}
\end{subequations}
with the matrices
\begin{subequations}
    \begin{align}
        \vec{M}_\rho           & = \int_{0}^{L} \vec{\Phi}\transp \rho A \vec{\Phi} \d s                 \eqcomma                  \\
        \vec{M}_S              & = \int_{0}^{L} \vec{\Psi} \transp \vec{\Psi} \d s                  \eqcomma                       \\
        \vec{F}_b              & = \int_{0}^{L} \vec{\Phi} \transp \vec{b} \d s                  \eqcomma                          \\
        \vec{K}(\hat{\vec{r}}) & = \int_{0}^{L} \vec{\Phi}_{,s}\transp \vec{\Phi}_{,s} \hat{\vec{r}}  \vec{\Psi} \d s       \eqdot
    \end{align}
\end{subequations}
Furthermore, $\vec{F}_N$ accounts for the forces due to Neumann boundary conditions in \eqref{weak_form_2}. Correspondingly, inserting the same approximations into the Hamiltonian \eqref{hamiltonian} yields a discretized version, which reads
\begin{align}
    H(\vec{r}\h, \vec{v}\h, C\h) & = \int_{0}^{L} \left( \frac{1}{2} \rho A \vec{v}\h \cdot \vec{v}\h + W(C\h) - \vec{r}\h \cdot \vec{b} \right) \d s                                          \notag \\
                                 & = \frac{1}{2} \hat{\vec{v}}\transp \vec{M}_\rho \hat{\vec{v}} + \int_{0}^{L} W(C\h) \d s - \hat{\vec{r}}\transp \vec{F}_b  \notag                                  \\
                                 & =: \hat{H}(\hat{\vec{r}},\hat{\vec{v}}, \hat{C}) \eqcomma \label{discrete_hamiltonian}
\end{align}
The discrete Hamiltonian $\hat{H}$ gives rise to the partial derivatives (i.e., discrete efforts corresponding to \eqref{efforts})
\begin{subequations}
    \begin{align}
        \frac{\partial \hat{H}}{\partial \hat{\vec{r}}} & = -\vec{F}_b \eqcomma                                   \\
        \frac{\partial \hat{H}}{\partial \hat{\vec{v}}} & = \vec{M}_\rho \hat{\vec{v}} \eqcomma                   \\
        \frac{\partial \hat{H}}{\partial \hat{\vec{C}}} & =  \int_{0}^{L} \vec{\Psi}\transp \D W(C\h) \d s \eqdot
    \end{align}
\end{subequations}
Taking into account \eqref{discrete_const1}, the last equation can be written as
\begin{align}
    \frac{\partial \hat{H}}{\partial \hat{\vec{C}}} = \vec{M}_S \frac{1}{2}\hat{\vec{S}} \eqdot \label{discrete_const_law}
\end{align}
This equation can be interpreted as discrete version of the constitutive relation \eqref{const_relation}.

We are now in the position to rewrite equations \eqref{eq:semi-discrete} governing the motion of the semi-discrete system in PH form. In particular, similar to \eqref{eq:PH-continuous-Mehrmann-Morandin2019}, we obtain
\begin{equation}
    \label{semidiscrete_EoM}
    \begin{array}{rcl}
        \vec{E}\frac{\d}{\d t}{\hat{\vec{x}}} & = & \vec{J} \hat{\vec{z}} + \vec{B}\hat{\vec{u}} \eqcomma \\
        \hat{\vec{y}}                         & = & \vec{B}\transp \hat{\vec{z}} \eqcomma
    \end{array}
    \qquad\mbox{where}\qquad
    \vec{E}\transp\hat{\vec{z}} = \D\hat{H}(\hat{\vec{x}}) \eqdot
\end{equation}
Or, in more detail,
\begin{align}
    \begin{bmatrix}
        \vec{I} & \vec{0}      & \vec{0}   \\
        \vec{0} & \vec{M}_\rho & \vec{0}   \\
        \vec{0} & \vec{0}      & \vec{M}_S
    \end{bmatrix}
    \frac{\d}{\d t}
    \begin{bmatrix}
        {\hat{\vec{r}}} \\ {\hat{\vec{v}}} \\ {\hat{\vec{C}}}
    \end{bmatrix} & =
    \begin{bmatrix}
        \vec{0}  &  & \vec{I}                          &  & \vec{0}                    \\
        -\vec{I} &  & \vec{0}                          &  & -2{\vec{K}}(\hat{\vec{r}}) \\
        \vec{0}  &  & 2{\vec{K}}(\hat{\vec{r}})\transp &  & \vec{0}
    \end{bmatrix}
    \begin{bmatrix}
        \partial_{\hat{\vec{r}}}\hat{H} \\
        \hat{\vec{v}}                   \\
        \frac{1}{2}\hat{\vec{S}}
    \end{bmatrix} +
    \begin{bmatrix}
        \vec{0} \\ \vec{B}_{\hat{\vec{v}}} \\ \vec{0}
    \end{bmatrix} \hat{\vec{u}} \eqcomma           \label{semidiscrete}                                 \\
    \hat{\vec{y}}                                                                                   & =
    \begin{bmatrix}
        \vec{0}\transp & \vec{B}_{\hat{\vec{v}}}\transp & \vec{0}\transp
    \end{bmatrix}
    \begin{bmatrix}
        \partial_{\hat{\vec{r}}}\hat{H} \\
        \hat{\vec{v}}                   \\
        \frac{1}{2}\hat{\vec{S}}
    \end{bmatrix} \eqdot
\end{align}
Here, $\vec{F}_{\mathrm{N}} = \vec{B}_{\hat{\vec{v}}} \hat{\vec{u}}$, where boundary inputs appear in the input vector $\hat{\vec{u}}$.
It can be observed that the matrix $\vec{J}$ is skew-symmetric ($\vec{J}\transp = - \vec{J}$) and $\vec{E}$ is symmetric ($\vec{E}\transp = \vec{E}$). The discrete power balance follows from
\begin{align}
    \frac{\d}{\d t}\hat{H}(\hat{\vec{x}})
    = \D\hat{H}(\hat{\vec{x}}) \cdot \frac{\d}{\d t}\hat{\vec{x}}
    = \vec{E}\transp\hat{\vec{z}} \cdot \frac{\d}{\d t}\hat{\vec{x}}
    = \hat{\vec{z}} \cdot \vec{E}\frac{\d}{\d t}\hat{\vec{x}}
    = \hat{\vec{z}} \cdot \left( \vec{J} \hat{\vec{z}} + \vec{B}\hat{\vec{u}} \right)
    = \hat{\vec{u}} \cdot \vec{B}\transp\hat{\vec{z}}
    = \hat{\vec{u}} \cdot \hat{\vec{y}} \eqdot \label{discrete_power}
\end{align}
Relation \eqref{discrete_power} expresses the passivity and losslessness of the spatially discrete system, which includes energy conservation for $\vec{u}=\vec{0}$.

\, \\
\noindent \textbf{Remark 4 (Discrete mixed boundary conditions):}
The application of mixed boundary conditions according to Remark 3 can be realized following standard finite element procedures. While Neumann boundary conditions are covered by $\vec{F}_{\mathrm{N}} = \vec{B}_{\hat{\vec{v}}} \hat{\vec{u}}$ in \eqref{semidiscrete}, Dirichlet boundary conditions can be enforced via Lagrange multipliers \cite{brugnoli_partitioned_2020} or in case of fixed bearings by directly setting the respective entries in the nodal vectors to the chosen values.

\subsection{Temporal discretization}
We aim at a structure-preserving time discretization of the semi-discrete PH system \eqref{semidiscrete_EoM}. To this end, we apply a one-step scheme which is closely related to the implicit mid-point rule. Let $\hat{\vec{x}}\n$ be the approximation of the state $\hat{\vec{x}}(t\n)$ at time $t\n$ and consider time steps of constant size $h = t\npe-t\n$. The time-stepping scheme
can now be written in the form
\begin{equation}
    \label{timediscrete_EoM}
    \begin{array}{rcl}
        \vec{E}\left(\hat{\vec{x}}\npe-\hat{\vec{x}}\n\right) & = & h \, \vec{J}(\hat{\vec{x}}\npeh) \hat{\vec{z}}\npeh + h \, \vec{B}\hat{\vec{u}}\npeh \eqcomma \\
        \hat{\vec{y}}\npeh                                    & = & \vec{B}\transp \hat{\vec{z}}\npeh \eqcomma
    \end{array}
\end{equation}
where $\hat{\vec{x}}\npeh=\frac{1}{2}(\hat{\vec{x}}\npe+\hat{\vec{x}}\n)$, $\hat{\vec{u}}\npeh\approx\hat{\vec{u}}(t\npeh)$, and $\hat{\vec{y}}\npeh\approx\hat{\vec{y}}(t\npeh)$. Moreover, $\hat{\vec{z}}\npeh$ is defined  through
\begin{equation}
    \vec{E}\transp\hat{\vec{z}}\npeh = \DG \hat{H}(\hat{\vec{x}}\n,\hat{\vec{x}}\npe) \eqcomma
\end{equation}
where $\DG \hat{H}(\hat{\vec{x}}\n,\hat{\vec{x}}\npe)$ is a discrete derivative in the sense of Gonzalez \cite{gonzalez_time_1996}. Among other properties (see \cite{gonzalez_time_1996} for more details), $\DG \hat{H}(\hat{\vec{x}}\n,\hat{\vec{x}}\npe)$ satisfies the crucial directionality property
\begin{equation} \label{eq:directionality}
    \DG \hat{H}(\hat{\vec{x}}\n,\hat{\vec{x}}\npe)\cdot\left(\hat{\vec{x}}\npe-\hat{\vec{x}}\n\right) = \hat{H}(\hat{\vec{x}}\npe) - \hat{H}(\hat{\vec{x}}\n) \eqdot
\end{equation}
In particular, we choose
\begin{equation} \label{discrete_derivative}
    \DG \hat{H}(\hat{\vec{x}}\n,\hat{\vec{x}}\npe) =
    \begin{bmatrix}
        \PDG_{\hat{\vec{r}}}\hat{H} \\
        \PDG_{\hat{\vec{v}}}\hat{H} \\
        \PDG_{\hat{\vec{C}}}\hat{H}
    \end{bmatrix} =
    \begin{bmatrix}
        -\vec{F}_b (t\npeh) \\
        \hat{\vec{v}}\npeh  \\
        \int_{0}^{L} \vec{\Psi} \transp \frac{W(C\h\npe) - W(C\h\n)}{C\h\npe - C\h\n} \d s
    \end{bmatrix} \eqcomma
\end{equation}
where the classical Greenspan's formula \cite{greenspan_conservative_1984} has been applied to arrive at a discrete derivative of $\D W (C\h)$. In the limit $C\h\npe \rightarrow C\h\n$, the mid-point evaluation of the standard derivative is used. We further note that the space integral in \eqref{discrete_derivative} is evaluated by means of Gaussian quadrature on element level using standard finite element assembly procedures.
It can be easily checked by a straightforward calculation that \eqref{discrete_derivative} does indeed satisfy the directionality property \eqref{eq:directionality} provided that $\vec{F}_b$ is constant.
Making use of the directionality property \eqref{eq:directionality} we obtain
\begin{align}
    \hat{H}\npe - \hat{H}\n \notag
    =\vec{E}\transp\hat{\vec{z}}\npeh\cdot\left(\hat{\vec{x}}\npe-\hat{\vec{x}}\n\right)
     & =\hat{\vec{z}}\npeh \cdot \vec{E}\left(\hat{\vec{x}}\npe-\hat{\vec{x}}\n\right)                                      \\
     & =\hat{\vec{z}}\npeh \cdot h \left( \vec{J}(\hat{\vec{x}}\npeh) \hat{\vec{z}}\npeh + \vec{B}\hat{\vec{u}}\npeh\right) \\
     & = h\, \hat{\vec{u}}\npeh \cdot \hat{\vec{y}}\npeh \eqcomma \notag
\end{align}
which is a time-discrete counterpart of \eqref{discrete_power}. This proves that the present time-stepping scheme exhibits passivity and losslessness (which includes energy-conservation in the case of vanishing inputs). We further remark without proof, that the present scheme also inherits symmetries of the underlying Hamiltonian, which implies conservation of the corresponding momentum maps.

\, \\
\noindent \textbf{Remark 5 (Discrete kinematic relation):}
The kinematic relation \eqref{rightCG} is exactly captured by the present method, i.e.,
\begin{align}
    C^{\mathrm{h}}\n = \partial_s \vec{r}^{\mathrm{h}}\n \cdot \partial_s \vec{r}^{\mathrm{h}}\n
\end{align}
for all $n$, if starting with consistent initial conditions satisfying $C^{\mathrm{h}}_0 = \partial_s \vec{r}^{\mathrm{h}}_0 \cdot \partial_s \vec{r}^{\mathrm{h}}_0$.

\section{Numerical example}
\label{sec_numerics}
We investigate the two-dimensional motion of a string made of rubber-like material ($E = 18400 \, \mathrm{N/m^2}$, $\rho=920 \, \mathrm{kg/m^3}$) with circular cross section (radius $R=0.0186 \, \mathrm{m}$). It is discretized in space with $n_{\mathrm{FE}}=30$ finite elements using linear ansatz functions for $\vec{\Phi}$ and constant ansatz functions for $\vec{\Psi}$. Numerical integration of the integrals has been achieved with two Gauss points per finite element and further simulation parameters are comprised in Table~\ref{tab:tab1}. We consider the initial conditions
\begin{subequations}
    \begin{align}
        \vec{r}_0(s) & = s \frac{\sqrt{2}}{2} \begin{bmatrix}
                                                  1 \\ -1
                                              \end{bmatrix} , \\
        \vec{v}_0(s) & = \vec{0}                ,             \\
        C_0(s)       & = 1 .
    \end{align}
\end{subequations}
The string is released at $t=0$ and simulated until $t=T$ with time step size $h$. In each time step \eqref{timediscrete_EoM} has been solved using Newton's method with a residual tolerance of $\epsilon_{\mathrm{Newton}}$. Considering mixed boundary conditions, the Dirichlet boundary conditions are given by fixing the left end of the string, such that
\begin{subequations}
    \begin{align}
        \vec{r}(s=0,t) & =\vec{0} , \\
        \vec{v}(s=0,t) & =\vec{0} ,
    \end{align}
\end{subequations}
and the Neumann boundary condition
\begin{align}
    \vec{F}_{\mathrm{N}}= \rho A  \begin{bmatrix}
                                      1 \\1
                                  \end{bmatrix} \sin \left( \pi \frac{t}{0.2}\right) \ , \qquad t \in \left[0 , 0.2 \right]
\end{align}
acts at the right end ($s=L$) of the string during a loading phase. After $t=0.2$, the system is closed and energy should be preserved during the motion. Throughout the simulation, gravitation is present such that
\begin{align}
    \vec{b} = -9.81 \rho A \begin{bmatrix}
                               0 \\ 1
                           \end{bmatrix} .
\end{align}
The material behavior is governed by the nonlinear relation \eqref{const_relation_hyper} and material constants simulating rubber, which can be found in Table~\ref{tab:tab1} as well.

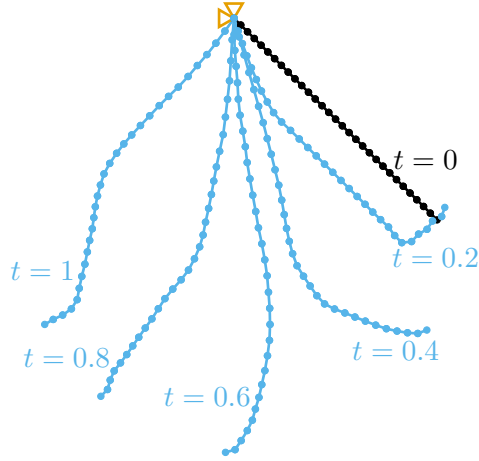
\begin{figure}[tb]
    \centering
    \setlength{\figH}{0.25\textheight}
    \setlength{\figW}{0.25\textheight}
    \input{images/conf_new.tikz}
    \caption{Initial configuration and snapshots during the motion}
    \label{fig:snapshots}
\end{figure}

\begin{table}[b]
    \begin{center}
        \caption{String pendulum example: Simulation parameters.}
        \vspace{1mm}
        \begin{tabular}{*{6}{|c}|}
            \hline
            $L \, [\mathrm{m}] $ & $EA\, [\mathrm{N}] $ & $\rho A \, [\mathrm{kg/m}] $ & $h \, [\mathrm{s}] $ & $T \, [\mathrm{s}] $ & $\epsilon_{\mathrm{Newton}} \, [-] $ \\
            \hline
            $1$                  & $20$                 & $1$                          & $1 \cdot 10^{-2}$    & $1$                  & $1 \cdot 10^{-11}$                   \\
            \hline
        \end{tabular}
        \label{tab:tab1}
    \end{center}
\end{table}

\begin{figure}[tb]

    \centering
    \setlength{\figH}{0.2\textheight}
    \setlength{\figW}{0.6\textwidth}
    \input{images/energy_new.tikz}
    \caption{Energy quantities}
    \label{fig:energy}

\end{figure}
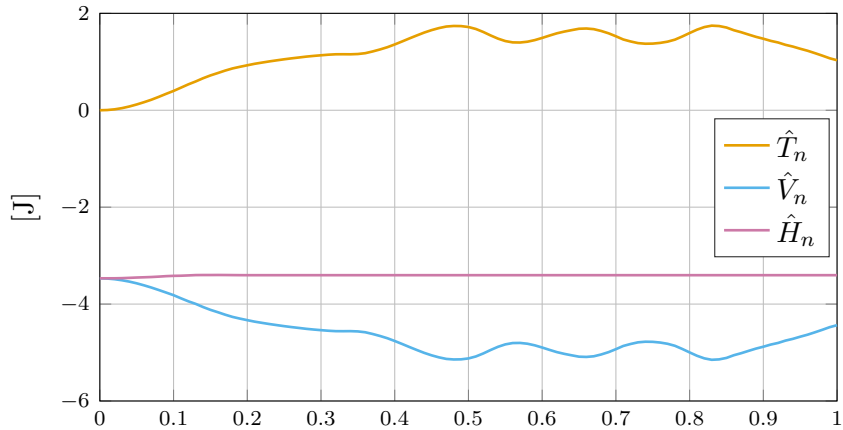

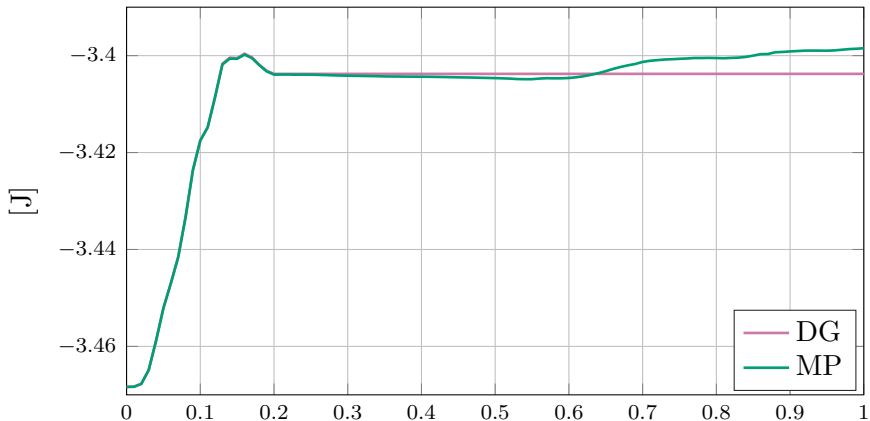
\begin{figure}[tb]
    \hspace*{-8mm}
    \centering
    \setlength{\figH}{0.2\textheight}
    \setlength{\figW}{0.6\textwidth}
    \input{images/H_comparison.tikz}
    \caption{Discrete Hamiltonian $\hat{H}\n$}
    \label{fig:Hamiltonian_comparison}

\end{figure}

\begin{figure}[tb]
    \centering

    \setlength{\figH}{0.2\textheight}
    \setlength{\figW}{0.6\textwidth}
    \input{images/Hdiff_new.tikz}
    \caption{Discrete Hamiltonian increments $| \hat{H}\npe - \hat{H}\n |$}
    \label{fig:energy_diff}
\end{figure}
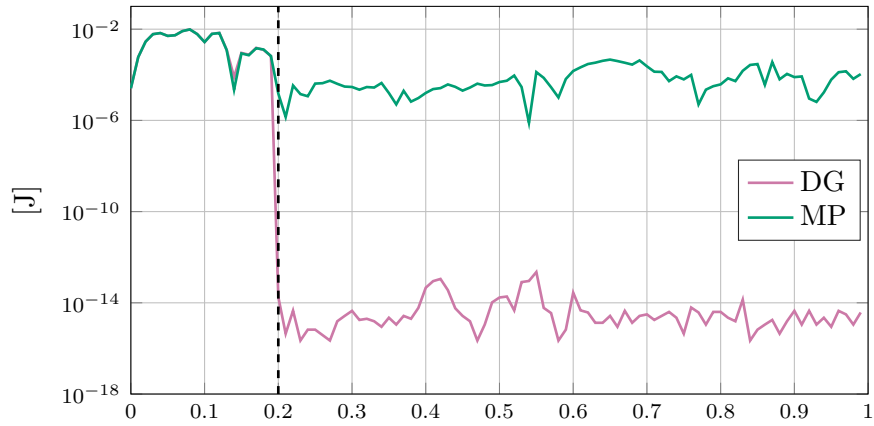

The evolution of the discrete Hamiltonian \eqref{discrete_hamiltonian} along with kinetic energy $\hat{T}\n$ and potential energy $\hat{V}\n$ (due to elastic deformations and gravity) are displayed in Fig.~\ref{fig:energy}. In Fig.~\ref{fig:Hamiltonian_comparison}, we compare our method (labeled \lq\lq DG\rq \rq) to the well-known implicit midpoint rule (labeled \lq\lq MP\rq \rq). During the loading phase, an increase of the Hamiltonian can be observed. This effect is accompanied by noticable time-step increments $| \hat{H}\npe - \hat{H}\n |$ for $t \in \left[0 , 0.2 \right]$, which are depicted in Fig.~\ref{fig:energy_diff}. After the loading phase, the discrete gradient approximation is able to capture the energy conservation numerically exactly, i.e., down to a level lower than the Newton tolerance. Once more, the results obtained with the midpoint discretization of \eqref{semidiscrete_EoM} are displayed for comparison (as used in \cite{thoma2022port}). Using the midpoint rule is not sufficient due to the nonlinear material \eqref{const_relation_hyper}, such that Hamiltonian increments are of the same order of magnitude as they are during loading.

\section{Conclusions}
\label{sec_conclusion}
In this work, a new port-Hamiltonian formulation for the analysis of geometrically exact strings with hyperelastic material behaviour has been proposed. The model features a nonlinear Hamiltonian formulated in position, velocity and strain quantities and therefore captures nonlinear material behaviour, extending previous work \cite{thoma2022port}. The typical skew-adjoint structure of infinite-dimensional PH systems is present in the governing state space equations. The passivity of the string formulation has been demonstrated with contact forces and velocities as boundary port variables. The weak form of the initial boundary value problem has been derived.

Based on the new formulation, we have derived a structure-preserving finite-dimensional PH state space model by using mixed finite elements. Therein, continuous ansatz functions for position and velocities and discontinuous ansatz functions for stress and strain quantities have been used to obtain the spatial discretization. This semi-discrete model inherits passivity (thus including energy-conservation) from the continuous string formulation and has subsequently been discretized in time in order to tackle simulation problems in an energy-consistent way. This property has been achieved due to the approximation with discrete gradients in the sense of Gonzalez \cite{gonzalez_time_1996}.

The numerical properties of the devised time-stepping scheme for nonlinear strings have been analyzed in a straightforward example problem. This has underlined the structure-preserving advantages of the devised method, in contrast to the well-known midpoint rule.

In future work, the newly proposed model might be used for model order reduction and control design. Besides the solution of the inverse dynamics (as, e.g., in \cite{strohle2022simultaneous}), state estimation and feedback control for the highly underactuated system will offer challenges for upcoming research activities.

%% file: images/filament.pdf_tex
\begingroup%
  \makeatletter%
  \providecommand\color[2][]{%
    \errmessage{(Inkscape) Color is used for the text in Inkscape, but the package 'color.sty' is not loaded}%
    \renewcommand\color[2][]{}%
  }%
  \providecommand\transparent[1]{%
    \errmessage{(Inkscape) Transparency is used (non-zero) for the text in Inkscape, but the package 'transparent.sty' is not loaded}%
    \renewcommand\transparent[1]{}%
  }%
  \providecommand\rotatebox[2]{#2}%
  \newcommand*\fsize{\dimexpr\f@size pt\relax}%
  \newcommand*\lineheight[1]{\fontsize{\fsize}{#1\fsize}\selectfont}%
  \ifx\svgwidth\undefined%
    \setlength{\unitlength}{841.88976378bp}%
    \ifx\svgscale\undefined%
      \relax%
    \else%
      \setlength{\unitlength}{\unitlength * \real{\svgscale}}%
    \fi%
  \else%
    \setlength{\unitlength}{\svgwidth}%
  \fi%
  \global\let\svgwidth\undefined%
  \global\let\svgscale\undefined%
  \makeatother%
  \begin{picture}(1,0.70707071)%
    \lineheight{1}%
    \setlength\tabcolsep{0pt}%
    \put(0,0){\includegraphics[width=\unitlength,page=1]{./images/filament.pdf}}%
    \put(0.043221,0.62557229){\color[rgb]{0,0,0}\makebox(0,0)[lt]{\lineheight{1.25}\smash{\begin{tabular}[t]{l}$X_2$\end{tabular}}}}%
    \put(0.24670978,0.19549103){\color[rgb]{0,0,0}\makebox(0,0)[lt]{\lineheight{1.25}\smash{\begin{tabular}[t]{l}$\mathbf{r}(s,t)$\end{tabular}}}}%
    \put(0.56273053,0.55874065){\color[rgb]{0,0,0}\makebox(0,0)[lt]{\lineheight{1.25}\smash{\begin{tabular}[t]{l}$\Omega_t$\end{tabular}}}}%
    \put(0.5984288,0.45833369){\color[rgb]{0,0,0}\makebox(0,0)[lt]{\lineheight{1.25}\smash{\begin{tabular}[t]{l}$\partial_s \mathbf{r}(s,t)$\end{tabular}}}}%
    \put(0.32854369,0.076994){\color[rgb]{0,0,0}\makebox(0,0)[lt]{\lineheight{1.25}\smash{\begin{tabular}[t]{l}$\Omega$\end{tabular}}}}%
    \put(0.85531971,0.06522028){\color[rgb]{0,0,0}\makebox(0,0)[lt]{\lineheight{1.25}\smash{\begin{tabular}[t]{l}$X_1$\end{tabular}}}}%
    \put(0,0){\includegraphics[width=\unitlength,page=2]{./images/filament.pdf}}%
    \put(0.5931638,0.08458091){\color[rgb]{0,0,0}\makebox(0,0)[lt]{\lineheight{1.25}\smash{\begin{tabular}[t]{l}$L$\end{tabular}}}}%
    \put(0,0){\includegraphics[width=\unitlength,page=3]{./images/filament.pdf}}%
  \end{picture}%
\endgroup%

%% file: images/conf_new.tikz
%
\begin{tikzpicture}

  \begin{axis}[%
      width=\figW,
      height=\figH,
      at={(0\figW,0\figH)},
      scale only axis,
      xmin=-0.85,
      xmax=0.85,
      ymin=-1.6,
      ymax=0.1,
      axis line style={draw=none},
      ticks=none,
      axis x line*=bottom,
      axis y line*=left
    ]

    \node[] at (0.6671,	-0.4971) {$t=0$};
    \node[] at (0.7016,	-0.8386) {\textcolor{color2}{$t=0.2$}};
    \node[] at (0.5633,	-1.167) {\textcolor{color2}{$t=0.4$}};
    \node[] at (-0.0931,	-1.326) {\textcolor{color2}{$t=0.6$}};
    \node[] at (-0.5821,	-1.19) {\textcolor{color2}{$t=0.8$}};
    \node[] at (-0.6661,	-0.885) {\textcolor{color2}{$t=1$}};

    \addplot[area legend, line width=1pt, table/row sep=crcr, patch, mesh, patch type=polygon, vertex count=2, color=black, fill=none, forget plot, patch table={%
          0	1\\
          2	3\\
          4	5\\
          6	7\\
          8	9\\
          10	11\\
          12	13\\
          14	15\\
          16	17\\
          18	19\\
          20	21\\
          22	23\\
          24	25\\
          26	27\\
          28	29\\
          30	31\\
          32	33\\
          34	35\\
          36	37\\
          38	39\\
          40	41\\
          42	43\\
          44	45\\
          46	47\\
          48	49\\
          50	51\\
          52	53\\
          54	55\\
          56	57\\
          58	59\\
        }]
    table[row sep=crcr] {%
        x	y\\
        0	-0\\
        0.02357023	-0.02357023\\
        0.02357023	-0.02357023\\
        0.04714045	-0.04714045\\
        0.04714045	-0.04714045\\
        0.07071068	-0.07071068\\
        0.07071068	-0.07071068\\
        0.0942809	-0.0942809\\
        0.0942809	-0.0942809\\
        0.1178511	-0.1178511\\
        0.1178511	-0.1178511\\
        0.1414214	-0.1414214\\
        0.1414214	-0.1414214\\
        0.1649916	-0.1649916\\
        0.1649916	-0.1649916\\
        0.1885618	-0.1885618\\
        0.1885618	-0.1885618\\
        0.212132	-0.212132\\
        0.212132	-0.212132\\
        0.2357023	-0.2357023\\
        0.2357023	-0.2357023\\
        0.2592725	-0.2592725\\
        0.2592725	-0.2592725\\
        0.2828427	-0.2828427\\
        0.2828427	-0.2828427\\
        0.3064129	-0.3064129\\
        0.3064129	-0.3064129\\
        0.3299832	-0.3299832\\
        0.3299832	-0.3299832\\
        0.3535534	-0.3535534\\
        0.3535534	-0.3535534\\
        0.3771236	-0.3771236\\
        0.3771236	-0.3771236\\
        0.4006938	-0.4006938\\
        0.4006938	-0.4006938\\
        0.4242641	-0.4242641\\
        0.4242641	-0.4242641\\
        0.4478343	-0.4478343\\
        0.4478343	-0.4478343\\
        0.4714045	-0.4714045\\
        0.4714045	-0.4714045\\
        0.4949747	-0.4949747\\
        0.4949747	-0.4949747\\
        0.518545	-0.518545\\
        0.518545	-0.518545\\
        0.5421152	-0.5421152\\
        0.5421152	-0.5421152\\
        0.5656854	-0.5656854\\
        0.5656854	-0.5656854\\
        0.5892557	-0.5892557\\
        0.5892557	-0.5892557\\
        0.6128259	-0.6128259\\
        0.6128259	-0.6128259\\
        0.6363961	-0.6363961\\
        0.6363961	-0.6363961\\
        0.6599663	-0.6599663\\
        0.6599663	-0.6599663\\
        0.6835366	-0.6835366\\
        0.6835366	-0.6835366\\
        0.7071068	-0.7071068\\
      };

    \addplot[area legend, mark size=1.2pt, mark=*, mark options={solid, black}, draw=none, table/row sep=crcr, patch, mesh, patch type=polygon, vertex count=2, color=black, fill=none, forget plot, patch table={%
          0	1\\
          2	3\\
          4	5\\
          6	7\\
          8	9\\
          10	11\\
          12	13\\
          14	15\\
          16	17\\
          18	19\\
          20	21\\
          22	23\\
          24	25\\
          26	27\\
          28	29\\
          30	31\\
          32	33\\
          34	35\\
          36	37\\
          38	39\\
          40	41\\
          42	43\\
          44	45\\
          46	47\\
          48	49\\
          50	51\\
          52	53\\
          54	55\\
          56	57\\
          58	59\\
        }]
    table[row sep=crcr] {%
        x	y\\
        0	-0\\
        0.02357023	-0.02357023\\
        0.02357023	-0.02357023\\
        0.04714045	-0.04714045\\
        0.04714045	-0.04714045\\
        0.07071068	-0.07071068\\
        0.07071068	-0.07071068\\
        0.0942809	-0.0942809\\
        0.0942809	-0.0942809\\
        0.1178511	-0.1178511\\
        0.1178511	-0.1178511\\
        0.1414214	-0.1414214\\
        0.1414214	-0.1414214\\
        0.1649916	-0.1649916\\
        0.1649916	-0.1649916\\
        0.1885618	-0.1885618\\
        0.1885618	-0.1885618\\
        0.212132	-0.212132\\
        0.212132	-0.212132\\
        0.2357023	-0.2357023\\
        0.2357023	-0.2357023\\
        0.2592725	-0.2592725\\
        0.2592725	-0.2592725\\
        0.2828427	-0.2828427\\
        0.2828427	-0.2828427\\
        0.3064129	-0.3064129\\
        0.3064129	-0.3064129\\
        0.3299832	-0.3299832\\
        0.3299832	-0.3299832\\
        0.3535534	-0.3535534\\
        0.3535534	-0.3535534\\
        0.3771236	-0.3771236\\
        0.3771236	-0.3771236\\
        0.4006938	-0.4006938\\
        0.4006938	-0.4006938\\
        0.4242641	-0.4242641\\
        0.4242641	-0.4242641\\
        0.4478343	-0.4478343\\
        0.4478343	-0.4478343\\
        0.4714045	-0.4714045\\
        0.4714045	-0.4714045\\
        0.4949747	-0.4949747\\
        0.4949747	-0.4949747\\
        0.518545	-0.518545\\
        0.518545	-0.518545\\
        0.5421152	-0.5421152\\
        0.5421152	-0.5421152\\
        0.5656854	-0.5656854\\
        0.5656854	-0.5656854\\
        0.5892557	-0.5892557\\
        0.5892557	-0.5892557\\
        0.6128259	-0.6128259\\
        0.6128259	-0.6128259\\
        0.6363961	-0.6363961\\
        0.6363961	-0.6363961\\
        0.6599663	-0.6599663\\
        0.6599663	-0.6599663\\
        0.6835366	-0.6835366\\
        0.6835366	-0.6835366\\
        0.7071068	-0.7071068\\
      };

    \addplot[area legend, line width=1pt, table/row sep=crcr, patch, mesh, patch type=polygon, vertex count=2, color=color2, fill=none, forget plot, patch table={%
          0	1\\
          2	3\\
          4	5\\
          6	7\\
          8	9\\
          10	11\\
          12	13\\
          14	15\\
          16	17\\
          18	19\\
          20	21\\
          22	23\\
          24	25\\
          26	27\\
          28	29\\
          30	31\\
          32	33\\
          34	35\\
          36	37\\
          38	39\\
          40	41\\
          42	43\\
          44	45\\
          46	47\\
          48	49\\
          50	51\\
          52	53\\
          54	55\\
          56	57\\
          58	59\\
        }]
    table[row sep=crcr] {%
        x	y\\
        0	0\\
        0.01437377	-0.04468353\\
        0.01437377	-0.04468353\\
        0.02916301	-0.08862133\\
        0.02916301	-0.08862133\\
        0.04465653	-0.1315046\\
        0.04465653	-0.1315046\\
        0.06053198	-0.1736198\\
        0.06053198	-0.1736198\\
        0.07775576	-0.2146684\\
        0.07775576	-0.2146684\\
        0.09545568	-0.2548471\\
        0.09545568	-0.2548471\\
        0.1158261	-0.2933656\\
        0.1158261	-0.2933656\\
        0.1377944	-0.3306752\\
        0.1377944	-0.3306752\\
        0.1656256	-0.3636027\\
        0.1656256	-0.3636027\\
        0.1968381	-0.3930875\\
        0.1968381	-0.3930875\\
        0.2271277	-0.4231646\\
        0.2271277	-0.4231646\\
        0.2568932	-0.4531433\\
        0.2568932	-0.4531433\\
        0.2862965	-0.4824895\\
        0.2862965	-0.4824895\\
        0.3152012	-0.5114006\\
        0.3152012	-0.5114006\\
        0.3438748	-0.5400765\\
        0.3438748	-0.5400765\\
        0.3724676	-0.5686628\\
        0.3724676	-0.5686628\\
        0.400832	-0.597043\\
        0.400832	-0.597043\\
        0.428744	-0.6249428\\
        0.428744	-0.6249428\\
        0.4562052	-0.6522973\\
        0.4562052	-0.6522973\\
        0.4830458	-0.6796411\\
        0.4830458	-0.6796411\\
        0.5103175	-0.7061988\\
        0.5103175	-0.7061988\\
        0.5379019	-0.7319948\\
        0.5379019	-0.7319948\\
        0.5593714	-0.7620435\\
        0.5593714	-0.7620435\\
        0.584311	-0.7870268\\
        0.584311	-0.7870268\\
        0.6182784	-0.7830967\\
        0.6182784	-0.7830967\\
        0.642694	-0.7568283\\
        0.642694	-0.7568283\\
        0.6720523	-0.7400699\\
        0.6720523	-0.7400699\\
        0.6913011	-0.7108405\\
        0.6913011	-0.7108405\\
        0.7229259	-0.6952384\\
        0.7229259	-0.6952384\\
        0.7339185	-0.6633033\\
      };

    \addplot[area legend, mark size=1.2pt, mark=*, mark options={solid, color2}, draw=none, table/row sep=crcr, patch, mesh, patch type=polygon, vertex count=2, color=color2, fill=none, forget plot, patch table={%
          0	1\\
          2	3\\
          4	5\\
          6	7\\
          8	9\\
          10	11\\
          12	13\\
          14	15\\
          16	17\\
          18	19\\
          20	21\\
          22	23\\
          24	25\\
          26	27\\
          28	29\\
          30	31\\
          32	33\\
          34	35\\
          36	37\\
          38	39\\
          40	41\\
          42	43\\
          44	45\\
          46	47\\
          48	49\\
          50	51\\
          52	53\\
          54	55\\
          56	57\\
          58	59\\
        }]
    table[row sep=crcr] {%
        x	y\\
        0	0\\
        0.01437377	-0.04468353\\
        0.01437377	-0.04468353\\
        0.02916301	-0.08862133\\
        0.02916301	-0.08862133\\
        0.04465653	-0.1315046\\
        0.04465653	-0.1315046\\
        0.06053198	-0.1736198\\
        0.06053198	-0.1736198\\
        0.07775576	-0.2146684\\
        0.07775576	-0.2146684\\
        0.09545568	-0.2548471\\
        0.09545568	-0.2548471\\
        0.1158261	-0.2933656\\
        0.1158261	-0.2933656\\
        0.1377944	-0.3306752\\
        0.1377944	-0.3306752\\
        0.1656256	-0.3636027\\
        0.1656256	-0.3636027\\
        0.1968381	-0.3930875\\
        0.1968381	-0.3930875\\
        0.2271277	-0.4231646\\
        0.2271277	-0.4231646\\
        0.2568932	-0.4531433\\
        0.2568932	-0.4531433\\
        0.2862965	-0.4824895\\
        0.2862965	-0.4824895\\
        0.3152012	-0.5114006\\
        0.3152012	-0.5114006\\
        0.3438748	-0.5400765\\
        0.3438748	-0.5400765\\
        0.3724676	-0.5686628\\
        0.3724676	-0.5686628\\
        0.400832	-0.597043\\
        0.400832	-0.597043\\
        0.428744	-0.6249428\\
        0.428744	-0.6249428\\
        0.4562052	-0.6522973\\
        0.4562052	-0.6522973\\
        0.4830458	-0.6796411\\
        0.4830458	-0.6796411\\
        0.5103175	-0.7061988\\
        0.5103175	-0.7061988\\
        0.5379019	-0.7319948\\
        0.5379019	-0.7319948\\
        0.5593714	-0.7620435\\
        0.5593714	-0.7620435\\
        0.584311	-0.7870268\\
        0.584311	-0.7870268\\
        0.6182784	-0.7830967\\
        0.6182784	-0.7830967\\
        0.642694	-0.7568283\\
        0.642694	-0.7568283\\
        0.6720523	-0.7400699\\
        0.6720523	-0.7400699\\
        0.6913011	-0.7108405\\
        0.6913011	-0.7108405\\
        0.7229259	-0.6952384\\
        0.7229259	-0.6952384\\
        0.7339185	-0.6633033\\
      };

    \addplot[area legend, line width=1pt, table/row sep=crcr, patch, mesh, patch type=polygon, vertex count=2, color=color2, fill=none, forget plot, patch table={%
          0	1\\
          2	3\\
          4	5\\
          6	7\\
          8	9\\
          10	11\\
          12	13\\
          14	15\\
          16	17\\
          18	19\\
          20	21\\
          22	23\\
          24	25\\
          26	27\\
          28	29\\
          30	31\\
          32	33\\
          34	35\\
          36	37\\
          38	39\\
          40	41\\
          42	43\\
          44	45\\
          46	47\\
          48	49\\
          50	51\\
          52	53\\
          54	55\\
          56	57\\
          58	59\\
        }]
    table[row sep=crcr] {%
        x	y\\
        0	0\\
        0.01779941	-0.06497641\\
        0.01779941	-0.06497641\\
        0.03528704	-0.1283052\\
        0.03528704	-0.1283052\\
        0.05218508	-0.1899502\\
        0.05218508	-0.1899502\\
        0.06924798	-0.2498864\\
        0.06924798	-0.2498864\\
        0.08564987	-0.3087319\\
        0.08564987	-0.3087319\\
        0.1015648	-0.3660227\\
        0.1015648	-0.3660227\\
        0.1156974	-0.4218656\\
        0.1156974	-0.4218656\\
        0.1295885	-0.4767097\\
        0.1295885	-0.4767097\\
        0.1428929	-0.5312251\\
        0.1428929	-0.5312251\\
        0.1558697	-0.5852564\\
        0.1558697	-0.5852564\\
        0.1661177	-0.6377364\\
        0.1661177	-0.6377364\\
        0.1748637	-0.6884481\\
        0.1748637	-0.6884481\\
        0.1829431	-0.7369959\\
        0.1829431	-0.7369959\\
        0.1925044	-0.7823845\\
        0.1925044	-0.7823845\\
        0.2040447	-0.8257573\\
        0.2040447	-0.8257573\\
        0.219447	-0.8678591\\
        0.219447	-0.8678591\\
        0.241983	-0.9045553\\
        0.241983	-0.9045553\\
        0.2649894	-0.9389482\\
        0.2649894	-0.9389482\\
        0.2885346	-0.9730059\\
        0.2885346	-0.9730059\\
        0.3153665	-1.00189\\
        0.3153665	-1.00189\\
        0.3500615	-1.019738\\
        0.3500615	-1.019738\\
        0.3859426	-1.033532\\
        0.3859426	-1.033532\\
        0.4211005	-1.047887\\
        0.4211005	-1.047887\\
        0.4572748	-1.062183\\
        0.4572748	-1.062183\\
        0.4934045	-1.076302\\
        0.4934045	-1.076302\\
        0.5304274	-1.088871\\
        0.5304274	-1.088871\\
        0.5675915	-1.097121\\
        0.5675915	-1.097121\\
        0.6040237	-1.10314\\
        0.6040237	-1.10314\\
        0.6397702	-1.105469\\
        0.6397702	-1.105469\\
        0.671312	-1.093531\\
      };

    \addplot[area legend, mark size=1.2pt, mark=*, mark options={solid, color2}, draw=none, table/row sep=crcr, patch, mesh, patch type=polygon, vertex count=2, color=color2, fill=none, forget plot, patch table={%
          0	1\\
          2	3\\
          4	5\\
          6	7\\
          8	9\\
          10	11\\
          12	13\\
          14	15\\
          16	17\\
          18	19\\
          20	21\\
          22	23\\
          24	25\\
          26	27\\
          28	29\\
          30	31\\
          32	33\\
          34	35\\
          36	37\\
          38	39\\
          40	41\\
          42	43\\
          44	45\\
          46	47\\
          48	49\\
          50	51\\
          52	53\\
          54	55\\
          56	57\\
          58	59\\
        }]
    table[row sep=crcr] {%
        x	y\\
        0	0\\
        0.01779941	-0.06497641\\
        0.01779941	-0.06497641\\
        0.03528704	-0.1283052\\
        0.03528704	-0.1283052\\
        0.05218508	-0.1899502\\
        0.05218508	-0.1899502\\
        0.06924798	-0.2498864\\
        0.06924798	-0.2498864\\
        0.08564987	-0.3087319\\
        0.08564987	-0.3087319\\
        0.1015648	-0.3660227\\
        0.1015648	-0.3660227\\
        0.1156974	-0.4218656\\
        0.1156974	-0.4218656\\
        0.1295885	-0.4767097\\
        0.1295885	-0.4767097\\
        0.1428929	-0.5312251\\
        0.1428929	-0.5312251\\
        0.1558697	-0.5852564\\
        0.1558697	-0.5852564\\
        0.1661177	-0.6377364\\
        0.1661177	-0.6377364\\
        0.1748637	-0.6884481\\
        0.1748637	-0.6884481\\
        0.1829431	-0.7369959\\
        0.1829431	-0.7369959\\
        0.1925044	-0.7823845\\
        0.1925044	-0.7823845\\
        0.2040447	-0.8257573\\
        0.2040447	-0.8257573\\
        0.219447	-0.8678591\\
        0.219447	-0.8678591\\
        0.241983	-0.9045553\\
        0.241983	-0.9045553\\
        0.2649894	-0.9389482\\
        0.2649894	-0.9389482\\
        0.2885346	-0.9730059\\
        0.2885346	-0.9730059\\
        0.3153665	-1.00189\\
        0.3153665	-1.00189\\
        0.3500615	-1.019738\\
        0.3500615	-1.019738\\
        0.3859426	-1.033532\\
        0.3859426	-1.033532\\
        0.4211005	-1.047887\\
        0.4211005	-1.047887\\
        0.4572748	-1.062183\\
        0.4572748	-1.062183\\
        0.4934045	-1.076302\\
        0.4934045	-1.076302\\
        0.5304274	-1.088871\\
        0.5304274	-1.088871\\
        0.5675915	-1.097121\\
        0.5675915	-1.097121\\
        0.6040237	-1.10314\\
        0.6040237	-1.10314\\
        0.6397702	-1.105469\\
        0.6397702	-1.105469\\
        0.671312	-1.093531\\
      };

    \addplot[area legend, line width=1pt, table/row sep=crcr, patch, mesh, patch type=polygon, vertex count=2, color=color2, fill=none, forget plot, patch table={%
          0	1\\
          2	3\\
          4	5\\
          6	7\\
          8	9\\
          10	11\\
          12	13\\
          14	15\\
          16	17\\
          18	19\\
          20	21\\
          22	23\\
          24	25\\
          26	27\\
          28	29\\
          30	31\\
          32	33\\
          34	35\\
          36	37\\
          38	39\\
          40	41\\
          42	43\\
          44	45\\
          46	47\\
          48	49\\
          50	51\\
          52	53\\
          54	55\\
          56	57\\
          58	59\\
        }]
    table[row sep=crcr] {%
        x	y\\
        0	0\\
        0.0003854613	-0.05699061\\
        0.0003854613	-0.05699061\\
        0.002922319	-0.1124393\\
        0.002922319	-0.1124393\\
        0.007339777	-0.1683551\\
        0.007339777	-0.1683551\\
        0.01217779	-0.2262211\\
        0.01217779	-0.2262211\\
        0.0150545	-0.2835128\\
        0.0150545	-0.2835128\\
        0.02112117	-0.3398084\\
        0.02112117	-0.3398084\\
        0.03068196	-0.3957295\\
        0.03068196	-0.3957295\\
        0.03935893	-0.4503478\\
        0.03935893	-0.4503478\\
        0.04589624	-0.5058723\\
        0.04589624	-0.5058723\\
        0.05530788	-0.5630216\\
        0.05530788	-0.5630216\\
        0.06519631	-0.6203623\\
        0.06519631	-0.6203623\\
        0.0755769	-0.6778119\\
        0.0755769	-0.6778119\\
        0.08505692	-0.7348531\\
        0.08505692	-0.7348531\\
        0.09696296	-0.7907527\\
        0.09696296	-0.7907527\\
        0.1071775	-0.8472308\\
        0.1071775	-0.8472308\\
        0.1153438	-0.9032065\\
        0.1153438	-0.9032065\\
        0.1199896	-0.960631\\
        0.1199896	-0.960631\\
        0.1242943	-1.01715\\
        0.1242943	-1.01715\\
        0.1254592	-1.073879\\
        0.1254592	-1.073879\\
        0.1234871	-1.129542\\
        0.1234871	-1.129542\\
        0.1165484	-1.18373\\
        0.1165484	-1.18373\\
        0.1075194	-1.235933\\
        0.1075194	-1.235933\\
        0.09673219	-1.284902\\
        0.09673219	-1.284902\\
        0.086885	-1.329789\\
        0.086885	-1.329789\\
        0.07443966	-1.370777\\
        0.07443966	-1.370777\\
        0.06390363	-1.409734\\
        0.06390363	-1.409734\\
        0.04629043	-1.446142\\
        0.04629043	-1.446142\\
        0.02487915	-1.481306\\
        0.02487915	-1.481306\\
        0.003717918	-1.513172\\
        0.003717918	-1.513172\\
        -0.02954671	-1.522672\\
      };

    \addplot[area legend, mark size=1.2pt, mark=*, mark options={solid, color2}, draw=none, table/row sep=crcr, patch, mesh, patch type=polygon, vertex count=2, color=color2, fill=none, forget plot, patch table={%
          0	1\\
          2	3\\
          4	5\\
          6	7\\
          8	9\\
          10	11\\
          12	13\\
          14	15\\
          16	17\\
          18	19\\
          20	21\\
          22	23\\
          24	25\\
          26	27\\
          28	29\\
          30	31\\
          32	33\\
          34	35\\
          36	37\\
          38	39\\
          40	41\\
          42	43\\
          44	45\\
          46	47\\
          48	49\\
          50	51\\
          52	53\\
          54	55\\
          56	57\\
          58	59\\
        }]
    table[row sep=crcr] {%
        x	y\\
        0	0\\
        0.0003854613	-0.05699061\\
        0.0003854613	-0.05699061\\
        0.002922319	-0.1124393\\
        0.002922319	-0.1124393\\
        0.007339777	-0.1683551\\
        0.007339777	-0.1683551\\
        0.01217779	-0.2262211\\
        0.01217779	-0.2262211\\
        0.0150545	-0.2835128\\
        0.0150545	-0.2835128\\
        0.02112117	-0.3398084\\
        0.02112117	-0.3398084\\
        0.03068196	-0.3957295\\
        0.03068196	-0.3957295\\
        0.03935893	-0.4503478\\
        0.03935893	-0.4503478\\
        0.04589624	-0.5058723\\
        0.04589624	-0.5058723\\
        0.05530788	-0.5630216\\
        0.05530788	-0.5630216\\
        0.06519631	-0.6203623\\
        0.06519631	-0.6203623\\
        0.0755769	-0.6778119\\
        0.0755769	-0.6778119\\
        0.08505692	-0.7348531\\
        0.08505692	-0.7348531\\
        0.09696296	-0.7907527\\
        0.09696296	-0.7907527\\
        0.1071775	-0.8472308\\
        0.1071775	-0.8472308\\
        0.1153438	-0.9032065\\
        0.1153438	-0.9032065\\
        0.1199896	-0.960631\\
        0.1199896	-0.960631\\
        0.1242943	-1.01715\\
        0.1242943	-1.01715\\
        0.1254592	-1.073879\\
        0.1254592	-1.073879\\
        0.1234871	-1.129542\\
        0.1234871	-1.129542\\
        0.1165484	-1.18373\\
        0.1165484	-1.18373\\
        0.1075194	-1.235933\\
        0.1075194	-1.235933\\
        0.09673219	-1.284902\\
        0.09673219	-1.284902\\
        0.086885	-1.329789\\
        0.086885	-1.329789\\
        0.07443966	-1.370777\\
        0.07443966	-1.370777\\
        0.06390363	-1.409734\\
        0.06390363	-1.409734\\
        0.04629043	-1.446142\\
        0.04629043	-1.446142\\
        0.02487915	-1.481306\\
        0.02487915	-1.481306\\
        0.003717918	-1.513172\\
        0.003717918	-1.513172\\
        -0.02954671	-1.522672\\
      };

    \addplot[area legend, line width=1pt, table/row sep=crcr, patch, mesh, patch type=polygon, vertex count=2, color=color2, fill=none, forget plot, patch table={%
          0	1\\
          2	3\\
          4	5\\
          6	7\\
          8	9\\
          10	11\\
          12	13\\
          14	15\\
          16	17\\
          18	19\\
          20	21\\
          22	23\\
          24	25\\
          26	27\\
          28	29\\
          30	31\\
          32	33\\
          34	35\\
          36	37\\
          38	39\\
          40	41\\
          42	43\\
          44	45\\
          46	47\\
          48	49\\
          50	51\\
          52	53\\
          54	55\\
          56	57\\
          58	59\\
        }]
    table[row sep=crcr] {%
        x	y\\
        0	0\\
        -0.006436379	-0.07527136\\
        -0.006436379	-0.07527136\\
        -0.01103297	-0.1490414\\
        -0.01103297	-0.1490414\\
        -0.01561639	-0.2191048\\
        -0.01561639	-0.2191048\\
        -0.022593	-0.2840087\\
        -0.022593	-0.2840087\\
        -0.02990546	-0.3464624\\
        -0.02990546	-0.3464624\\
        -0.03769122	-0.4079683\\
        -0.03769122	-0.4079683\\
        -0.05009069	-0.4662063\\
        -0.05009069	-0.4662063\\
        -0.06129387	-0.5211966\\
        -0.06129387	-0.5211966\\
        -0.07047681	-0.5707517\\
        -0.07047681	-0.5707517\\
        -0.07864774	-0.6179901\\
        -0.07864774	-0.6179901\\
        -0.08831242	-0.6668051\\
        -0.08831242	-0.6668051\\
        -0.09909351	-0.7149682\\
        -0.09909351	-0.7149682\\
        -0.1076027	-0.7600837\\
        -0.1076027	-0.7600837\\
        -0.1195579	-0.8015968\\
        -0.1195579	-0.8015968\\
        -0.1376286	-0.8435402\\
        -0.1376286	-0.8435402\\
        -0.1584856	-0.8830406\\
        -0.1584856	-0.8830406\\
        -0.1841278	-0.9172085\\
        -0.1841278	-0.9172085\\
        -0.2131304	-0.9497847\\
        -0.2131304	-0.9497847\\
        -0.2405211	-0.9833714\\
        -0.2405211	-0.9833714\\
        -0.2635069	-1.018732\\
        -0.2635069	-1.018732\\
        -0.2858313	-1.05425\\
        -0.2858313	-1.05425\\
        -0.3089821	-1.086471\\
        -0.3089821	-1.086471\\
        -0.3315941	-1.118729\\
        -0.3315941	-1.118729\\
        -0.3525117	-1.149126\\
        -0.3525117	-1.149126\\
        -0.373313	-1.179539\\
        -0.373313	-1.179539\\
        -0.3952568	-1.208801\\
        -0.3952568	-1.208801\\
        -0.4171949	-1.236854\\
        -0.4171949	-1.236854\\
        -0.4310732	-1.269145\\
        -0.4310732	-1.269145\\
        -0.4404317	-1.302141\\
        -0.4404317	-1.302141\\
        -0.4635511	-1.325977\\
      };

    \addplot[area legend, mark size=1.2pt, mark=*, mark options={solid, color2}, draw=none, table/row sep=crcr, patch, mesh, patch type=polygon, vertex count=2, color=color2, fill=none, forget plot, patch table={%
          0	1\\
          2	3\\
          4	5\\
          6	7\\
          8	9\\
          10	11\\
          12	13\\
          14	15\\
          16	17\\
          18	19\\
          20	21\\
          22	23\\
          24	25\\
          26	27\\
          28	29\\
          30	31\\
          32	33\\
          34	35\\
          36	37\\
          38	39\\
          40	41\\
          42	43\\
          44	45\\
          46	47\\
          48	49\\
          50	51\\
          52	53\\
          54	55\\
          56	57\\
          58	59\\
        }]
    table[row sep=crcr] {%
        x	y\\
        0	0\\
        -0.006436379	-0.07527136\\
        -0.006436379	-0.07527136\\
        -0.01103297	-0.1490414\\
        -0.01103297	-0.1490414\\
        -0.01561639	-0.2191048\\
        -0.01561639	-0.2191048\\
        -0.022593	-0.2840087\\
        -0.022593	-0.2840087\\
        -0.02990546	-0.3464624\\
        -0.02990546	-0.3464624\\
        -0.03769122	-0.4079683\\
        -0.03769122	-0.4079683\\
        -0.05009069	-0.4662063\\
        -0.05009069	-0.4662063\\
        -0.06129387	-0.5211966\\
        -0.06129387	-0.5211966\\
        -0.07047681	-0.5707517\\
        -0.07047681	-0.5707517\\
        -0.07864774	-0.6179901\\
        -0.07864774	-0.6179901\\
        -0.08831242	-0.6668051\\
        -0.08831242	-0.6668051\\
        -0.09909351	-0.7149682\\
        -0.09909351	-0.7149682\\
        -0.1076027	-0.7600837\\
        -0.1076027	-0.7600837\\
        -0.1195579	-0.8015968\\
        -0.1195579	-0.8015968\\
        -0.1376286	-0.8435402\\
        -0.1376286	-0.8435402\\
        -0.1584856	-0.8830406\\
        -0.1584856	-0.8830406\\
        -0.1841278	-0.9172085\\
        -0.1841278	-0.9172085\\
        -0.2131304	-0.9497847\\
        -0.2131304	-0.9497847\\
        -0.2405211	-0.9833714\\
        -0.2405211	-0.9833714\\
        -0.2635069	-1.018732\\
        -0.2635069	-1.018732\\
        -0.2858313	-1.05425\\
        -0.2858313	-1.05425\\
        -0.3089821	-1.086471\\
        -0.3089821	-1.086471\\
        -0.3315941	-1.118729\\
        -0.3315941	-1.118729\\
        -0.3525117	-1.149126\\
        -0.3525117	-1.149126\\
        -0.373313	-1.179539\\
        -0.373313	-1.179539\\
        -0.3952568	-1.208801\\
        -0.3952568	-1.208801\\
        -0.4171949	-1.236854\\
        -0.4171949	-1.236854\\
        -0.4310732	-1.269145\\
        -0.4310732	-1.269145\\
        -0.4404317	-1.302141\\
        -0.4404317	-1.302141\\
        -0.4635511	-1.325977\\
      };

    \addplot[area legend, line width=1pt, table/row sep=crcr, patch, mesh, patch type=polygon, vertex count=2, color=color2, fill=none, forget plot, patch table={%
          0	1\\
          2	3\\
          4	5\\
          6	7\\
          8	9\\
          10	11\\
          12	13\\
          14	15\\
          16	17\\
          18	19\\
          20	21\\
          22	23\\
          24	25\\
          26	27\\
          28	29\\
          30	31\\
          32	33\\
          34	35\\
          36	37\\
          38	39\\
          40	41\\
          42	43\\
          44	45\\
          46	47\\
          48	49\\
          50	51\\
          52	53\\
          54	55\\
          56	57\\
          58	59\\
        }]
    table[row sep=crcr] {%
        x	y\\
        0	0\\
        -0.0315147	-0.0447051\\
        -0.0315147	-0.0447051\\
        -0.06303566	-0.08713187\\
        -0.06303566	-0.08713187\\
        -0.09422648	-0.1271506\\
        -0.09422648	-0.1271506\\
        -0.1276191	-0.1655256\\
        -0.1276191	-0.1655256\\
        -0.159793	-0.2039264\\
        -0.159793	-0.2039264\\
        -0.1939657	-0.2410368\\
        -0.1939657	-0.2410368\\
        -0.2277905	-0.2738055\\
        -0.2277905	-0.2738055\\
        -0.2603162	-0.306188\\
        -0.2603162	-0.306188\\
        -0.2934961	-0.3419041\\
        -0.2934961	-0.3419041\\
        -0.3257475	-0.3759772\\
        -0.3257475	-0.3759772\\
        -0.3561054	-0.4090597\\
        -0.3561054	-0.4090597\\
        -0.3836394	-0.4405365\\
        -0.3836394	-0.4405365\\
        -0.4085996	-0.4732693\\
        -0.4085996	-0.4732693\\
        -0.4326335	-0.5088714\\
        -0.4326335	-0.5088714\\
        -0.4526486	-0.5457498\\
        -0.4526486	-0.5457498\\
        -0.4650342	-0.5848016\\
        -0.4650342	-0.5848016\\
        -0.4753031	-0.6233408\\
        -0.4753031	-0.6233408\\
        -0.4828493	-0.6609763\\
        -0.4828493	-0.6609763\\
        -0.4880435	-0.6992076\\
        -0.4880435	-0.6992076\\
        -0.4934526	-0.7382222\\
        -0.4934526	-0.7382222\\
        -0.5025268	-0.7788931\\
        -0.5025268	-0.7788931\\
        -0.5106175	-0.8195311\\
        -0.5106175	-0.8195311\\
        -0.5191362	-0.8625298\\
        -0.5191362	-0.8625298\\
        -0.5299103	-0.9053365\\
        -0.5299103	-0.9053365\\
        -0.538833	-0.9471301\\
        -0.538833	-0.9471301\\
        -0.5461194	-0.9870375\\
        -0.5461194	-0.9870375\\
        -0.5661586	-1.019288\\
        -0.5661586	-1.019288\\
        -0.5965126	-1.040897\\
        -0.5965126	-1.040897\\
        -0.6293675	-1.056331\\
        -0.6293675	-1.056331\\
        -0.6590096	-1.073855\\
      };

    \addplot[area legend, mark size=1.2pt, mark=*, mark options={solid, color2}, draw=none, table/row sep=crcr, patch, mesh, patch type=polygon, vertex count=2, color=color2, fill=none, forget plot, patch table={%
          0	1\\
          2	3\\
          4	5\\
          6	7\\
          8	9\\
          10	11\\
          12	13\\
          14	15\\
          16	17\\
          18	19\\
          20	21\\
          22	23\\
          24	25\\
          26	27\\
          28	29\\
          30	31\\
          32	33\\
          34	35\\
          36	37\\
          38	39\\
          40	41\\
          42	43\\
          44	45\\
          46	47\\
          48	49\\
          50	51\\
          52	53\\
          54	55\\
          56	57\\
          58	59\\
        }]
    table[row sep=crcr] {%
        x	y\\
        0	0\\
        -0.0315147	-0.0447051\\
        -0.0315147	-0.0447051\\
        -0.06303566	-0.08713187\\
        -0.06303566	-0.08713187\\
        -0.09422648	-0.1271506\\
        -0.09422648	-0.1271506\\
        -0.1276191	-0.1655256\\
        -0.1276191	-0.1655256\\
        -0.159793	-0.2039264\\
        -0.159793	-0.2039264\\
        -0.1939657	-0.2410368\\
        -0.1939657	-0.2410368\\
        -0.2277905	-0.2738055\\
        -0.2277905	-0.2738055\\
        -0.2603162	-0.306188\\
        -0.2603162	-0.306188\\
        -0.2934961	-0.3419041\\
        -0.2934961	-0.3419041\\
        -0.3257475	-0.3759772\\
        -0.3257475	-0.3759772\\
        -0.3561054	-0.4090597\\
        -0.3561054	-0.4090597\\
        -0.3836394	-0.4405365\\
        -0.3836394	-0.4405365\\
        -0.4085996	-0.4732693\\
        -0.4085996	-0.4732693\\
        -0.4326335	-0.5088714\\
        -0.4326335	-0.5088714\\
        -0.4526486	-0.5457498\\
        -0.4526486	-0.5457498\\
        -0.4650342	-0.5848016\\
        -0.4650342	-0.5848016\\
        -0.4753031	-0.6233408\\
        -0.4753031	-0.6233408\\
        -0.4828493	-0.6609763\\
        -0.4828493	-0.6609763\\
        -0.4880435	-0.6992076\\
        -0.4880435	-0.6992076\\
        -0.4934526	-0.7382222\\
        -0.4934526	-0.7382222\\
        -0.5025268	-0.7788931\\
        -0.5025268	-0.7788931\\
        -0.5106175	-0.8195311\\
        -0.5106175	-0.8195311\\
        -0.5191362	-0.8625298\\
        -0.5191362	-0.8625298\\
        -0.5299103	-0.9053365\\
        -0.5299103	-0.9053365\\
        -0.538833	-0.9471301\\
        -0.538833	-0.9471301\\
        -0.5461194	-0.9870375\\
        -0.5461194	-0.9870375\\
        -0.5661586	-1.019288\\
        -0.5661586	-1.019288\\
        -0.5965126	-1.040897\\
        -0.5965126	-1.040897\\
        -0.6293675	-1.056331\\
        -0.6293675	-1.056331\\
        -0.6590096	-1.073855\\
      };
    \addplot [color=color1, line width=1pt, forget plot]
    table[row sep=crcr]{%
        0	0\\
        -0.05196	0.03\\
        -0.05196	-0.03\\
        0	0\\
      };
    \addplot [color=color1, line width=1pt, forget plot]
    table[row sep=crcr]{%
        0	0\\
        -0.03	0.05196\\
        0.03	0.05196\\
        0	0\\
      };
  \end{axis}
\end{tikzpicture}%

%% file: images/energy_new.tikz
%
\definecolor{mycolor1}{rgb}{0.00000,0.44700,0.74100}%
\definecolor{mycolor2}{rgb}{0.85000,0.32500,0.09800}%
\definecolor{mycolor3}{rgb}{0.92900,0.69400,0.12500}%
\begin{tikzpicture}

  \begin{axis}[%
      width=0.951\figW,
      height=\figH,
      at={(0\figW,0\figH)},
      scale only axis,
      xmin=0,
      xmax=1,
      ymin=-6,
      ymax=2,
      xmajorgrids,
      ymajorgrids,
      ylabel={[J]},
      axis background/.style={fill=white},
      legend style={at={(0.99,0.55)}, anchor = east, legend cell align=left, align=left, draw=white!15!black}
    ]
    \addplot [color=color1, line width=1pt]
    table[row sep=crcr]{%
        0	0\\
        0.01	0.00469161\\
        0.02	0.01892293\\
        0.03	0.04349219\\
        0.04	0.0775423\\
        0.05	0.1194501\\
        0.06	0.1665173\\
        0.07	0.2177588\\
        0.08	0.2768101\\
        0.09	0.3391396\\
        0.1	0.3989017\\
        0.11	0.4660352\\
        0.12	0.534111\\
        0.13	0.5942102\\
        0.14	0.6605574\\
        0.15	0.7184468\\
        0.16	0.7678959\\
        0.17	0.8186376\\
        0.18	0.8629338\\
        0.19	0.8967276\\
        0.2	0.928054\\
        0.21	0.9587977\\
        0.22	0.9838052\\
        0.23	1.008263\\
        0.24	1.030198\\
        0.25	1.051655\\
        0.26	1.071077\\
        0.27	1.088708\\
        0.28	1.106882\\
        0.29	1.121607\\
        0.3	1.134436\\
        0.31	1.146617\\
        0.32	1.155016\\
        0.33	1.155887\\
        0.34	1.15478\\
        0.35	1.160387\\
        0.36	1.177699\\
        0.37	1.213207\\
        0.38	1.254422\\
        0.39	1.301119\\
        0.4	1.358125\\
        0.41	1.419412\\
        0.42	1.482818\\
        0.43	1.543964\\
        0.44	1.603198\\
        0.45	1.65382\\
        0.46	1.695865\\
        0.47	1.727824\\
        0.48	1.739279\\
        0.49	1.736792\\
        0.5	1.717785\\
        0.51	1.677654\\
        0.52	1.616064\\
        0.53	1.544567\\
        0.54	1.478466\\
        0.55	1.428022\\
        0.56	1.400556\\
        0.57	1.397971\\
        0.58	1.414442\\
        0.59	1.447797\\
        0.6	1.492673\\
        0.61	1.541552\\
        0.62	1.587963\\
        0.63	1.623501\\
        0.64	1.653473\\
        0.65	1.681754\\
        0.66	1.686869\\
        0.67	1.6737\\
        0.68	1.636478\\
        0.69	1.5884\\
        0.7	1.526975\\
        0.71	1.460601\\
        0.72	1.419446\\
        0.73	1.385746\\
        0.74	1.372853\\
        0.75	1.375486\\
        0.76	1.390015\\
        0.77	1.415305\\
        0.78	1.456806\\
        0.79	1.522087\\
        0.8	1.59372\\
        0.81	1.662924\\
        0.82	1.716687\\
        0.83	1.746011\\
        0.84	1.737769\\
        0.85	1.707319\\
        0.86	1.65152\\
        0.87	1.609838\\
        0.88	1.561755\\
        0.89	1.514989\\
        0.9	1.47618\\
        0.91	1.432029\\
        0.92	1.398345\\
        0.93	1.350152\\
        0.94	1.314259\\
        0.95	1.269462\\
        0.96	1.22796\\
        0.97	1.178435\\
        0.98	1.127072\\
        0.99	1.072691\\
        1	1.031913\\
      };
    \addlegendentry{$\hat{T}\n$}

    \addplot [color=color2, line width=1pt]
    table[row sep=crcr]{%
        0	-3.468359\\
        0.01	-3.473025\\
        0.02	-3.486657\\
        0.03	-3.508379\\
        0.04	-3.536374\\
        0.05	-3.571533\\
        0.06	-3.613524\\
        0.07	-3.659338\\
        0.08	-3.710153\\
        0.09	-3.762719\\
        0.1	-3.816378\\
        0.11	-3.880806\\
        0.12	-3.942661\\
        0.13	-3.995876\\
        0.14	-4.060978\\
        0.15	-4.118947\\
        0.16	-4.167485\\
        0.17	-4.218987\\
        0.18	-4.264802\\
        0.19	-4.299827\\
        0.2	-4.331816\\
        0.21	-4.362559\\
        0.22	-4.387567\\
        0.23	-4.412025\\
        0.24	-4.433959\\
        0.25	-4.455416\\
        0.26	-4.474839\\
        0.27	-4.49247\\
        0.28	-4.510644\\
        0.29	-4.525368\\
        0.3	-4.538198\\
        0.31	-4.550379\\
        0.32	-4.558778\\
        0.33	-4.559649\\
        0.34	-4.558541\\
        0.35	-4.564148\\
        0.36	-4.58146\\
        0.37	-4.616969\\
        0.38	-4.658183\\
        0.39	-4.70488\\
        0.4	-4.761886\\
        0.41	-4.823174\\
        0.42	-4.88658\\
        0.43	-4.947725\\
        0.44	-5.00696\\
        0.45	-5.057582\\
        0.46	-5.099626\\
        0.47	-5.131586\\
        0.48	-5.14304\\
        0.49	-5.140554\\
        0.5	-5.121547\\
        0.51	-5.081416\\
        0.52	-5.019825\\
        0.53	-4.948328\\
        0.54	-4.882228\\
        0.55	-4.831784\\
        0.56	-4.804318\\
        0.57	-4.801732\\
        0.58	-4.818204\\
        0.59	-4.851559\\
        0.6	-4.896435\\
        0.61	-4.945313\\
        0.62	-4.991725\\
        0.63	-5.027263\\
        0.64	-5.057234\\
        0.65	-5.085516\\
        0.66	-5.09063\\
        0.67	-5.077462\\
        0.68	-5.040239\\
        0.69	-4.992161\\
        0.7	-4.930737\\
        0.71	-4.864363\\
        0.72	-4.823208\\
        0.73	-4.789507\\
        0.74	-4.776615\\
        0.75	-4.779247\\
        0.76	-4.793777\\
        0.77	-4.819067\\
        0.78	-4.860568\\
        0.79	-4.925849\\
        0.8	-4.997482\\
        0.81	-5.066685\\
        0.82	-5.120448\\
        0.83	-5.149773\\
        0.84	-5.14153\\
        0.85	-5.111081\\
        0.86	-5.055282\\
        0.87	-5.013599\\
        0.88	-4.965517\\
        0.89	-4.918751\\
        0.9	-4.879942\\
        0.91	-4.835791\\
        0.92	-4.802106\\
        0.93	-4.753914\\
        0.94	-4.718021\\
        0.95	-4.673224\\
        0.96	-4.631721\\
        0.97	-4.582197\\
        0.98	-4.530833\\
        0.99	-4.476452\\
        1	-4.435675\\
      };
    \addlegendentry{$\hat{V}\n$}

    \addplot [color=color3, line width=1pt]
    table[row sep=crcr]{%
        0	-3.468359\\
        0.01	-3.468333\\
        0.02	-3.467734\\
        0.03	-3.464886\\
        0.04	-3.458832\\
        0.05	-3.452083\\
        0.06	-3.447007\\
        0.07	-3.441579\\
        0.08	-3.433343\\
        0.09	-3.423579\\
        0.1	-3.417477\\
        0.11	-3.41477\\
        0.12	-3.40855\\
        0.13	-3.401666\\
        0.14	-3.400421\\
        0.15	-3.4005\\
        0.16	-3.399589\\
        0.17	-3.40035\\
        0.18	-3.401868\\
        0.19	-3.403099\\
        0.2	-3.403762\\
        0.21	-3.403762\\
        0.22	-3.403762\\
        0.23	-3.403762\\
        0.24	-3.403762\\
        0.25	-3.403762\\
        0.26	-3.403762\\
        0.27	-3.403762\\
        0.28	-3.403762\\
        0.29	-3.403762\\
        0.3	-3.403762\\
        0.31	-3.403762\\
        0.32	-3.403762\\
        0.33	-3.403762\\
        0.34	-3.403762\\
        0.35	-3.403762\\
        0.36	-3.403762\\
        0.37	-3.403762\\
        0.38	-3.403762\\
        0.39	-3.403762\\
        0.4	-3.403762\\
        0.41	-3.403762\\
        0.42	-3.403762\\
        0.43	-3.403762\\
        0.44	-3.403762\\
        0.45	-3.403762\\
        0.46	-3.403762\\
        0.47	-3.403762\\
        0.48	-3.403762\\
        0.49	-3.403762\\
        0.5	-3.403762\\
        0.51	-3.403762\\
        0.52	-3.403762\\
        0.53	-3.403762\\
        0.54	-3.403762\\
        0.55	-3.403762\\
        0.56	-3.403762\\
        0.57	-3.403762\\
        0.58	-3.403762\\
        0.59	-3.403762\\
        0.6	-3.403762\\
        0.61	-3.403762\\
        0.62	-3.403762\\
        0.63	-3.403762\\
        0.64	-3.403762\\
        0.65	-3.403762\\
        0.66	-3.403762\\
        0.67	-3.403762\\
        0.68	-3.403762\\
        0.69	-3.403762\\
        0.7	-3.403762\\
        0.71	-3.403762\\
        0.72	-3.403762\\
        0.73	-3.403762\\
        0.74	-3.403762\\
        0.75	-3.403762\\
        0.76	-3.403762\\
        0.77	-3.403762\\
        0.78	-3.403762\\
        0.79	-3.403762\\
        0.8	-3.403762\\
        0.81	-3.403762\\
        0.82	-3.403762\\
        0.83	-3.403762\\
        0.84	-3.403762\\
        0.85	-3.403762\\
        0.86	-3.403762\\
        0.87	-3.403762\\
        0.88	-3.403762\\
        0.89	-3.403762\\
        0.9	-3.403762\\
        0.91	-3.403762\\
        0.92	-3.403762\\
        0.93	-3.403762\\
        0.94	-3.403762\\
        0.95	-3.403762\\
        0.96	-3.403762\\
        0.97	-3.403762\\
        0.98	-3.403762\\
        0.99	-3.403762\\
        1	-3.403762\\
      };
    \addlegendentry{$\hat{H}\n$}

  \end{axis}
\end{tikzpicture}%

%% file: images/H_comparison.tikz
%
%
\begin{tikzpicture}

  \begin{axis}[%
      width=0.951\figW,
      height=\figH,
      at={(0\figW,0\figH)},
      scale only axis,
      xmin=0,
      xmax=1,
      ymin=-3.47,
      ymax=-3.39,
      xmajorgrids,
      ymajorgrids,
      ylabel={[J]},
      axis background/.style={fill=white},
      legend style={at={(0.99,0.01)}, anchor = south east, legend cell align=left, align=left, draw=white!15!black}
    ]

    \addplot [color=color3, line width=1pt]
    table[row sep=crcr]{%
        0	-3.468359\\
        0.01	-3.468333\\
        0.02	-3.467734\\
        0.03	-3.464886\\
        0.04	-3.458832\\
        0.05	-3.452083\\
        0.06	-3.447007\\
        0.07	-3.441579\\
        0.08	-3.433343\\
        0.09	-3.423579\\
        0.1	-3.417477\\
        0.11	-3.41477\\
        0.12	-3.40855\\
        0.13	-3.401666\\
        0.14	-3.400421\\
        0.15	-3.4005\\
        0.16	-3.399589\\
        0.17	-3.40035\\
        0.18	-3.401868\\
        0.19	-3.403099\\
        0.2	-3.403762\\
        0.21	-3.403762\\
        0.22	-3.403762\\
        0.23	-3.403762\\
        0.24	-3.403762\\
        0.25	-3.403762\\
        0.26	-3.403762\\
        0.27	-3.403762\\
        0.28	-3.403762\\
        0.29	-3.403762\\
        0.3	-3.403762\\
        0.31	-3.403762\\
        0.32	-3.403762\\
        0.33	-3.403762\\
        0.34	-3.403762\\
        0.35	-3.403762\\
        0.36	-3.403762\\
        0.37	-3.403762\\
        0.38	-3.403762\\
        0.39	-3.403762\\
        0.4	-3.403762\\
        0.41	-3.403762\\
        0.42	-3.403762\\
        0.43	-3.403762\\
        0.44	-3.403762\\
        0.45	-3.403762\\
        0.46	-3.403762\\
        0.47	-3.403762\\
        0.48	-3.403762\\
        0.49	-3.403762\\
        0.5	-3.403762\\
        0.51	-3.403762\\
        0.52	-3.403762\\
        0.53	-3.403762\\
        0.54	-3.403762\\
        0.55	-3.403762\\
        0.56	-3.403762\\
        0.57	-3.403762\\
        0.58	-3.403762\\
        0.59	-3.403762\\
        0.6	-3.403762\\
        0.61	-3.403762\\
        0.62	-3.403762\\
        0.63	-3.403762\\
        0.64	-3.403762\\
        0.65	-3.403762\\
        0.66	-3.403762\\
        0.67	-3.403762\\
        0.68	-3.403762\\
        0.69	-3.403762\\
        0.7	-3.403762\\
        0.71	-3.403762\\
        0.72	-3.403762\\
        0.73	-3.403762\\
        0.74	-3.403762\\
        0.75	-3.403762\\
        0.76	-3.403762\\
        0.77	-3.403762\\
        0.78	-3.403762\\
        0.79	-3.403762\\
        0.8	-3.403762\\
        0.81	-3.403762\\
        0.82	-3.403762\\
        0.83	-3.403762\\
        0.84	-3.403762\\
        0.85	-3.403762\\
        0.86	-3.403762\\
        0.87	-3.403762\\
        0.88	-3.403762\\
        0.89	-3.403762\\
        0.9	-3.403762\\
        0.91	-3.403762\\
        0.92	-3.403762\\
        0.93	-3.403762\\
        0.94	-3.403762\\
        0.95	-3.403762\\
        0.96	-3.403762\\
        0.97	-3.403762\\
        0.98	-3.403762\\
        0.99	-3.403762\\
        1	-3.403762\\
      };
    \addlegendentry{DG}

    \addplot [color=color4, line width=1pt]
    table[row sep=crcr]{%
        0	-3.468359\\
        0.01	-3.468333\\
        0.02	-3.467735\\
        0.03	-3.46489\\
        0.04	-3.458846\\
        0.05	-3.452097\\
        0.06	-3.447024\\
        0.07	-3.441613\\
        0.08	-3.433375\\
        0.09	-3.423626\\
        0.1	-3.417546\\
        0.11	-3.414821\\
        0.12	-3.408597\\
        0.13	-3.401829\\
        0.14	-3.400626\\
        0.15	-3.400648\\
        0.16	-3.399786\\
        0.17	-3.400505\\
        0.18	-3.401964\\
        0.19	-3.40325\\
        0.2	-3.403912\\
        0.21	-3.403897\\
        0.22	-3.403895\\
        0.23	-3.403929\\
        0.24	-3.403915\\
        0.25	-3.403926\\
        0.26	-3.403968\\
        0.27	-3.40401\\
        0.28	-3.404065\\
        0.29	-3.404105\\
        0.3	-3.404136\\
        0.31	-3.404165\\
        0.32	-3.404187\\
        0.33	-3.404216\\
        0.34	-3.404244\\
        0.35	-3.404287\\
        0.36	-3.404303\\
        0.37	-3.404308\\
        0.38	-3.404328\\
        0.39	-3.404335\\
        0.4	-3.404344\\
        0.41	-3.40436\\
        0.42	-3.404384\\
        0.43	-3.40441\\
        0.44	-3.404447\\
        0.45	-3.404477\\
        0.46	-3.404497\\
        0.47	-3.404525\\
        0.48	-3.404566\\
        0.49	-3.4046\\
        0.5	-3.404635\\
        0.51	-3.404683\\
        0.52	-3.404738\\
        0.53	-3.404831\\
        0.54	-3.40486\\
        0.55	-3.40486\\
        0.56	-3.404729\\
        0.57	-3.404656\\
        0.58	-3.404684\\
        0.59	-3.404674\\
        0.6	-3.404609\\
        0.61	-3.404462\\
        0.62	-3.40425\\
        0.63	-3.403954\\
        0.64	-3.403614\\
        0.65	-3.403204\\
        0.66	-3.402741\\
        0.67	-3.402341\\
        0.68	-3.402003\\
        0.69	-3.401724\\
        0.7	-3.401297\\
        0.71	-3.401062\\
        0.72	-3.400926\\
        0.73	-3.400792\\
        0.74	-3.40074\\
        0.75	-3.400655\\
        0.76	-3.400592\\
        0.77	-3.400493\\
        0.78	-3.400488\\
        0.79	-3.400466\\
        0.8	-3.400497\\
        0.81	-3.400535\\
        0.82	-3.400465\\
        0.83	-3.400413\\
        0.84	-3.400263\\
        0.85	-3.399993\\
        0.86	-3.399698\\
        0.87	-3.399662\\
        0.88	-3.399302\\
        0.89	-3.399239\\
        0.9	-3.39913\\
        0.91	-3.399051\\
        0.92	-3.398966\\
        0.93	-3.398957\\
        0.94	-3.398963\\
        0.95	-3.39898\\
        0.96	-3.398919\\
        0.97	-3.39879\\
        0.98	-3.398648\\
        0.99	-3.398581\\
        1	-3.398474\\
      };
    \addlegendentry{MP}

  \end{axis}
\end{tikzpicture}%

%% file: images/Hdiff_new.tikz
%
\definecolor{mycolor1}{rgb}{0.00000,0.44700,0.74100}%
\begin{tikzpicture}

  \begin{axis}[%
      width=0.951\figW,
      height=\figH,
      at={(0\figW,0\figH)},
      scale only axis,
      xmin=0,
      xmax=1,
      ymin=1e-18,
      ymax=0.1,
      ymode=log,
      xmajorgrids,
      ymajorgrids,
      ylabel={[J]},
      yminorticks=true,
      axis background/.style={fill=white},
      legend style={at={(0.99,0.5)}, anchor = east, legend cell align=left, align=left, draw=white!15!black}
    ]
    \addplot [color=color3, line width=1pt]
    table[row sep=crcr]{%
        0	2.548e-05\\
        0.01	0.0005991\\
        0.02	0.002848\\
        0.03	0.006055\\
        0.04	0.006748\\
        0.05	0.005076\\
        0.06	0.005428\\
        0.07	0.008236\\
        0.08	0.009763\\
        0.09	0.006103\\
        0.1	0.002706\\
        0.11	0.00622\\
        0.12	0.006884\\
        0.13	0.001245\\
        0.14	7.929e-05\\
        0.15	0.0009116\\
        0.16	0.000761\\
        0.17	0.001519\\
        0.18	0.001231\\
        0.19	0.0006621\\
        0.2	1.799e-14\\
        0.21	4.441e-16\\
        0.22	4.552e-15\\
        0.23	2.22e-16\\
        0.24	6.661e-16\\
        0.25	6.661e-16\\
        0.26	0\\
        0.27	2.22e-16\\
        0.28	1.554e-15\\
        0.29	2.665e-15\\
        0.3	4.441e-15\\
        0.31	1.776e-15\\
        0.32	1.998e-15\\
        0.33	1.554e-15\\
        0.34	8.882e-16\\
        0.35	2.22e-15\\
        0.36	1.11e-15\\
        0.37	2.665e-15\\
        0.38	1.998e-15\\
        0.39	5.995e-15\\
        0.4	4.596e-14\\
        0.41	8.926e-14\\
        0.42	1.115e-13\\
        0.43	3.597e-14\\
        0.44	5.773e-15\\
        0.45	2.665e-15\\
        0.46	1.554e-15\\
        0.47	2.22e-16\\
        0.48	1.11e-15\\
        0.49	1.044e-14\\
        0.5	1.71e-14\\
        0.51	1.887e-14\\
        0.52	4.663e-15\\
        0.53	8.105e-14\\
        0.54	9.082e-14\\
        0.55	2.238e-13\\
        0.56	5.995e-15\\
        0.57	3.553e-15\\
        0.58	2.22e-16\\
        0.59	6.661e-16\\
        0.6	2.798e-14\\
        0.61	4.663e-15\\
        0.62	3.775e-15\\
        0.63	1.332e-15\\
        0.64	1.332e-15\\
        0.65	2.665e-15\\
        0.66	8.882e-16\\
        0.67	4.441e-15\\
        0.68	1.332e-15\\
        0.69	2.665e-15\\
        0.7	3.109e-15\\
        0.71	1.776e-15\\
        0.72	2.665e-15\\
        0.73	3.997e-15\\
        0.74	2.22e-15\\
        0.75	4.441e-16\\
        0.76	6.217e-15\\
        0.77	3.775e-15\\
        0.78	1.11e-15\\
        0.79	3.997e-15\\
        0.8	3.997e-15\\
        0.81	2.22e-15\\
        0.82	1.554e-15\\
        0.83	1.377e-14\\
        0.84	2.22e-16\\
        0.85	6.661e-16\\
        0.86	1.11e-15\\
        0.87	1.776e-15\\
        0.88	4.441e-16\\
        0.89	1.554e-15\\
        0.9	4.441e-15\\
        0.91	1.11e-15\\
        0.92	4.441e-15\\
        0.93	1.11e-15\\
        0.94	2.22e-15\\
        0.95	8.882e-16\\
        0.96	4.441e-15\\
        0.97	3.109e-15\\
        0.98	1.11e-15\\
        0.99	3.775e-15\\
      };
    \addlegendentry{DG}

    \addplot [color=color4, line width=1pt]
    table[row sep=crcr]{%
        0	2.541e-05\\
        0.01	0.0005984\\
        0.02	0.002845\\
        0.03	0.006044\\
        0.04	0.006749\\
        0.05	0.005073\\
        0.06	0.005412\\
        0.07	0.008238\\
        0.08	0.009749\\
        0.09	0.006081\\
        0.1	0.002725\\
        0.11	0.006224\\
        0.12	0.006768\\
        0.13	0.001203\\
        0.14	2.202e-05\\
        0.15	0.0008625\\
        0.16	0.0007194\\
        0.17	0.001459\\
        0.18	0.001286\\
        0.19	0.0006616\\
        0.2	1.503e-05\\
        0.21	1.409e-06\\
        0.22	3.383e-05\\
        0.23	1.42e-05\\
        0.24	1.145e-05\\
        0.25	4.109e-05\\
        0.26	4.273e-05\\
        0.27	5.494e-05\\
        0.28	4.019e-05\\
        0.29	3.039e-05\\
        0.3	2.918e-05\\
        0.31	2.216e-05\\
        0.32	2.902e-05\\
        0.33	2.758e-05\\
        0.34	4.359e-05\\
        0.35	1.584e-05\\
        0.36	5.023e-06\\
        0.37	1.979e-05\\
        0.38	6.637e-06\\
        0.39	9.542e-06\\
        0.4	1.596e-05\\
        0.41	2.36e-05\\
        0.42	2.599e-05\\
        0.43	3.773e-05\\
        0.44	2.964e-05\\
        0.45	2.017e-05\\
        0.46	2.753e-05\\
        0.47	4.078e-05\\
        0.48	3.397e-05\\
        0.49	3.535e-05\\
        0.5	4.809e-05\\
        0.51	5.487e-05\\
        0.52	9.265e-05\\
        0.53	2.905e-05\\
        0.54	7.913e-07\\
        0.55	0.0001311\\
        0.56	7.336e-05\\
        0.57	2.84e-05\\
        0.58	1.009e-05\\
        0.59	6.54e-05\\
        0.6	0.0001466\\
        0.61	0.0002123\\
        0.62	0.0002961\\
        0.63	0.0003395\\
        0.64	0.0004101\\
        0.65	0.0004634\\
        0.66	0.0003998\\
        0.67	0.0003376\\
        0.68	0.0002795\\
        0.69	0.0004267\\
        0.7	0.0002351\\
        0.71	0.0001364\\
        0.72	0.0001335\\
        0.73	5.263e-05\\
        0.74	8.521e-05\\
        0.75	6.281e-05\\
        0.76	9.865e-05\\
        0.77	5.014e-06\\
        0.78	2.226e-05\\
        0.79	3.157e-05\\
        0.8	3.789e-05\\
        0.81	6.986e-05\\
        0.82	5.271e-05\\
        0.83	0.0001496\\
        0.84	0.00027\\
        0.85	0.0002951\\
        0.86	3.624e-05\\
        0.87	0.0003595\\
        0.88	6.343e-05\\
        0.89	0.0001093\\
        0.9	7.892e-05\\
        0.91	8.462e-05\\
        0.92	9.056e-06\\
        0.93	6.461e-06\\
        0.94	1.694e-05\\
        0.95	6.101e-05\\
        0.96	0.0001295\\
        0.97	0.0001422\\
        0.98	6.652e-05\\
        0.99	0.0001071\\
      };
    \addlegendentry{MP}
    \addplot [color=black, dashed, line width=1pt, forget plot]
    table[row sep=crcr]{%
        0.2	1e-18\\
        0.2	1e-1\\
      };

  \end{axis}
\end{tikzpicture}%